\documentclass[11pt]{article}
\usepackage{amsmath}
\usepackage{amssymb}
\usepackage[margin=1in]{geometry}
\usepackage{multirow}
\usepackage{siunitx}
\usepackage{comment}

\usepackage{graphicx} 
\usepackage{authblk} 

\newcommand\bm[1]{\begin{bmatrix}#1\end{bmatrix}} 
\usepackage{amsthm}

\title{Arbitrary High Order Low-rank Completely Positive and Trace
Preserving (CPTP) Schemes for Lindblad Equations with Time-dependent Hamiltonian}
\author[1]{Jiuhua Hu\thanks{Department of Mathematics, Virginia Tech, Blacksburg, VA 24061 U.S.A. (\textsf{jiuhuahu@vt.edu}).}}
\author[1]{Daniel Appel\"{o}\thanks{Department of Mathematics, Virginia Tech, Blacksburg, VA 24061 U.S.A. (\textsf{appelo@vt.edu}). Research is supported by the U.S. Department of Energy, Office of Science, Advanced Scientific Computing Research (ASCR), under Award Number DE-SC0025424. This material is based upon work supported by the National Science Foundation under Grant No.
DMS-2436319}}
\author[1]{Yingda Cheng\thanks{Department of Mathematics, Virginia Tech, Blacksburg, VA 24061 U.S.A. (\textsf{yingda@vt.edu}). Research is supported by DOE grant DE-SC0023164, AFOSR grant FA9550-25-1-0154 and Virginia Tech. This material is based upon work supported by the National Science Foundation under Grant No.
DMS-2424139 while the second and third authors were in residence at the Simons Laufer Mathematical Sciences
Institute in Berkeley, California, during the Fall 2025 semester.}}
\date{\today}
\affil[1]{Department of Mathematics, Virginia Tech, Blacksburg, USA}

\begin{document}
\maketitle
\section*{Abstract}
In this paper, we develop a framework for designing arbitrary high order  low-rank schemes for  the  Lindblad equation with time-dependent Hamiltonians.
Our approach is based on nested Picard iterative integrators (NPI) and results in schemes in Kraus form that are completely positive and trace preserving (CPTP). The schemes are amenable to low rank formulations, making them suitable for problems where the matrix rank of the density matrix is small.      

\section{Introduction}
We are interested in numerical simulations of the dynamics of an open quantum system, described by the Lindblad equation \cite{lindblad1976generators}.  In particular, we consider systems with general time-dependent Hamiltonians which are important, e.g. for quantum device simulation and optimal control \cite{HOHO}.

The Lindblad equation evolves the density matrix, $\rho(t)\in \mathbb{C}^{d\times d}$,   according to  
\begin{align}\label{eqn_master}
    \frac{d \rho(t)}{d t} &= - i (H(t)\rho(t)-\rho(t)H(t)) + \mathcal{L}\rho(t),\nonumber\\
    \rho(0)&= \rho_0.
\end{align}
Here the initial state of the system is given by the density matrix $\rho_0$, and $H(t)= H_d+ H_c(t)$ is the system Hamiltonian consisting of a time-independent part $H_d$ describing the inherent properties of the quantum system and a time-dependent part $H_c(t)$ used to control the system. The Lindbladian operator $\mathcal{L} $ models system-environment interactions and takes the form:
\begin{equation}
    \mathcal{L}\rho(t) = \displaystyle \sum_{\alpha} \mathcal{L}_{\alpha} \rho(t) \mathcal{L}^{\dagger}_{\alpha}  - \frac{1}{2} \left( \mathcal{L}^{\dagger}_{\alpha}\mathcal{L}_{\alpha} \rho(t) + \rho(t) \mathcal{L}^{\dagger}_{\alpha}\mathcal{L}_{\alpha}
 \right).
\end{equation}
Here $\{\mathcal{L}_{\alpha}\}_{\alpha}  $ represent a set of jump operators that describe the dissipative part of the dynamics, and $ \dagger $ denotes the adjoint (conjugate transpose) of an operator. The density matrix $\rho(t)$ is positive semi-definite with unit trace. 

A defining feature of the Lindblad equation is that its evolution of $\rho$ preserves complete positivity (CP) and is trace preserving (TP) \cite{manzano2020short}. A goal of this paper is to construct numerical methods that preserve these properties. Choi's theorem (see e.g. \cite{manzano2020short}) states: {\em A linear map $\mathcal{G}$ is CP iff it can be represented as}
$$
\mathcal{G}\rho=\sum_l G_l^\dagger \rho G_l,
$$
and this provides a way to construct numerical methods that preserve CP. Let $\rho^n \approx \rho(t_n)$ be an approximation to the density matrix at time $t_n$ and let a numerical method be defined as a {\em Kraus map} $\mathcal{G}$ that takes $\rho^n$ to $\rho^{n+1}$, that is 
\[
\rho^{n+1} \equiv \mathcal{G}\rho^n=\sum_l G_l^\dagger \rho^n G_l.
\] 
Such a  method is obviously CP and preservation of trace can be achieved by either a  non-linear trace re-normalization \mbox{$\rho^{n+1} \leftarrow \rho^{n+1}/{\rm Tr}(\rho^{n+1})$} or the recently proposed linear trace renormalization procedure from  \cite{robin2025unconditionally}. 

When quantum systems (for example multiple qubits) are coupled,  the dimension of the total quantum state space is equal to the product of the dimension of each subsystem. Thus, in many quantum computing applications the dimension of the quantum state space $d$ can become very large. Managing such large state spaces is challenging due to the need to store and operate on the density matrix. Addressing this issue requires exploiting dimension reduction techniques. Here we are principally interested in the case when the open system at hand is close to a closed system. This is the case for quantum computing applications. When the system is closed,  the Lindbladian operator $\mathcal{L}\rho = 0$. Then the Lindblad equation \eqref{eqn_master} is reduced to the von Neumann equation whose solution has rank one.  

In our prior work \cite{appelo2024kraus}, we developed high-order low-rank CPTP schemes for Lindblad equation with time-independent Hamiltonians. The methods there  belong to the class of integrating factor methods (or Lawson methods \cite{lawson1967generalized}). For the schemes in \cite{appelo2024kraus} the CP property is met as long as the entries in Butcher tableaus are non-negative. In \cite{appelo2024kraus} we use the truncated SVD, which we also prove is a CP map, for low-rank compression. The methods in \cite{appelo2024kraus} enjoy simplicity in design and provable CPTP property. However, they do not work for time-dependent Hamiltonians.

The goal of this paper is to design new high-order CPTP schemes for Lindblad equations for time-dependent Hamiltonians with CPTP property. Most schemes in the numerical literature satisfying CPTP property follows a reformulation of the equation introduced in \cite{steinbach1995high}. This includes \cite{cao2025structure, chen2024full,appelo2024kraus,chen2025full}, among which, \cite{chen2024full, chen2025full} developed exponential Euler and exponential midpoint methods, which are applicable to time-dependent Hamiltonian cases. Our formulation is based on  the Duhamel's principle from the reformulation of Lindblad equation in \cite{steinbach1995high}. This is similar to the approaches in \cite{cao2025structure, chen2024full,appelo2024kraus,chen2025full}, however, here we employ the nested Picard iterative integrator (NPI) framework \cite{baumstark2018uniformly,schratz2021low,cai2019uniformly}, which can easily be adopted to construct CPTP schemes of arbitrary high order of accuracy.  Low-rank compression of the density matrix is obtained by the use of the truncated SVD. We remark that the low-rank NPI schemes developed in this paper are closely connected to a well-known technique in numerical ODE: spectral deferred correction (SDC) \cite{dutt2000spectral}, which was recently explored in the low-rank context in \cite{li2024high,bachmayr2025iterative}.  

Lindblad dynamics as a \emph{linear} CPTP map possess many interesting physical properties, including the contractiveness of the nuclear norm of the density matrix and the von-Neumann entropy. Moreover, the eigenvalues of the Lindblad superoperator all have non-positive real part. In principle, if the discrete numerical operator is also a {\em linear} CPTP map, those properties can automatically be retained. However, since our method is not a \emph{linear} CPTP map, it necessitates a study of the stability property of our method. To understand the stability properties of the proposed methods, we conduct a stability analysis for single-qubit model problems. In this analysis we find that for a system that involves a single jump operator that models decay, the stability condition is determined only by the highest order flow operator. When the model includes dephasing jump operators, the stability condition is more involved but the stability properties are approximately the same as for the decay case when the dephasing is weak.

The rest of the paper is organized as follows. In Section \ref{sec: scheme}, we describe the proposed NPI schemes.  In Section \ref{sec:stability}, we carry out stability analysis on the family of methods for single-qubit model problems. Section \ref{sec:lowrank} briefly describes the low-rank truncation. In Section \ref{sec: numeric}, we present numerical examples including: coupled qudits and cavities with time dependent Hamiltonians, reference problems where an analytic solution is known, and a more applied example where we explore control function design in the Jaynes–Cummings model to actively suppress quantum revivals. We conclude the paper in Section \ref{sec: conclusions}.

\section{ CPTP Schemes by NPI} \label{sec: scheme}
In this section, we  formulate our high-order low-rank  CPTP schemes. Following \cite{steinbach1995high}, we  rewrite \eqref{eqn_master} into 
\begin{equation}
   \frac{d \rho(t)}{d t} = J(t) \rho(t) +\rho(t) J^\dagger(t)+ \sum_{\alpha} \mathcal{L}_{\alpha} \rho(t) \mathcal{L}^{\dagger}_{\alpha}, \label{eqn: master_rewrite}
\end{equation}
where $J(t) = -i H_{\text{eff}}$, and  $H_{\text{eff}}=H(t)+\frac{1}{2i}  \sum_{\alpha} \mathcal{L}^{\dagger}_{\alpha}\mathcal{L}_{\alpha}$. To simplify the notation we  introduce $\mathcal{L}_L \rho (t) = \sum_{\alpha} \mathcal{L}_{\alpha} \rho(t) \mathcal{L}^{\dagger}_{\alpha} $. 

Noting that the $\mathcal{L}_L \rho (t)$ term is already on Kraus form we momentarily and for purpose of illustration set it to zero. Thus we consider 
\begin{equation}
    \frac{d \rho(t)}{d t} = J(t) \rho(t) +\rho(t) J^\dagger(t).\label{model: rho_homogene}
\end{equation}
Since the right hand side of this equation is not on Kraus form it is not possible to directly use a standard ODE solver (e.g. an explicit Runge-Kutta method) and still have the CPTP property. However, as the positive semi-definite $\rho(t)$  admits a  decomposition  $\rho(t) = V(t) V^\dagger(t)$ solving for $\rho(t)$ satisfying equation \eqref{model: rho_homogene} with initial condition $\rho(t^0)=\rho_0$ is equivalent to solving for a matrix $V(t)$ satisfying
\begin{align}\label{eqn: homogeneous}
    \frac{dV(t)}{dt} &= J(t) V(t),
    \\
    V(t^0)&=V_0,\nonumber
\end{align}
where $V_0 $ is the matrix satisfying $\rho_0 = V_0 V_0^\dagger$.

Formally we can express the solution of \eqref{eqn: homogeneous} as $V(t) = U(t,s) V(s)$ where we refer to $U(t,s)$ as the {\em flow operator} evolving from time $s$ to time $t$. Consequently, $\rho(t)$ (solving (\ref{eqn: homogeneous})) can be written using the flow operator as
\begin{equation}
    \rho(t) =  U(t,s) \rho(s) U^\dagger(t,s). \label{eqn: rho}
\end{equation}

This observation motivates the use of Duhamel's principle to formally solve equation \eqref{eqn: master_rewrite} (which does contain $\mathcal{L}_L \rho (t)$). This allows for the representation of $\rho(t)$ evolving from an initial condition at time $\tau$ as 
\begin{eqnarray*}
\rho(t) = U(t,\tau) \rho(\tau) U^\dagger(t,\tau) + \int_\tau^t U(t,s) \mathcal{L}_L \rho(s) U^\dagger(t,s) \, ds \\
= U(t,\tau) \rho(\tau) U^\dagger(t,\tau) + \sum_{\alpha} \int_\tau^t [U(t,s)  \mathcal{L}_{\alpha} ] \rho(s)  [U(t,s) \mathcal{L}_{\alpha}]^\dagger \, ds.
\end{eqnarray*}

\subsection{Achieving Arbitrary-order of Accuracy by NPI}
\label{sec: NPI}

Let the time interval $ t \in [0,T]$ be partitioned by a uniform mesh $t^n = n \Delta t,\,  n = 0,\ldots,N_T,$ $\Delta t = T/N_T$. Then, using the above formula, the solution at time $t^{n+1}$, evolving from time $t^n$, is expressed as
\begin{equation}\label{eqn: rho_n+1}
    \rho(t^{n+1}) = U(t^{n+1},t^n) \rho(t^n) U^\dagger(t^{n+1},t^n) + \int_{t^n}^{t^{n+1}} U(t^{n+1},s) \mathcal{L}_L \rho(s) U^\dagger(t^{n+1},s) \, ds.
\end{equation}

Given a solution at time $t^n$ obtaining an approximation to the solution $\rho(t^{n+1})$ at a later time relies on: i.) approximating the flow operator $U(t^{n+1},t^n)$ by a numerical ODE solver and, ii.) approximating the integral by quadrature. These approximations will be discussed later but we note here that the Kraus representation structure of  $\rho(t^{n+1} )$ in \eqref{eqn: rho_n+1} will still be present for the discretized problem if the quadrature has positive weights. Note also that the Kraus form is independent of the approximation of the flow operator.

If the expression inside the integral (the integrand) was known for all times, and could be explicitly evaluated, constructing a CPTP scheme would amount to: i.) approximating the flow applied to the initial data $\rho(t^n)$ in the first term, and ii.) approximating the integral by a quadrature rule with positive weights. However, since the integrand is not known such a method would become implicit and solving for the solution at the next timestep would destroy the Kraus form and the method would, in general, not have the CPTP property. The ideas in \cite{cao2025structure} relies on repeated use of Duhamel's principle in \eqref{eqn: rho_n+1} which result in a single stage method. Here we present a simpler and more general implementation that uses NPI. The resulting method is of multi-stage type and iteratively increases the order of accuracy of the numerical solution by one order per iteration.

We now define our method to find an approximation of arbitrary order of accuracy to the solution $\rho(t^{n+1})$ expressed in the form \eqref{eqn: rho_n+1} and starting from a given $\rho^n$. The method will output a sequence of numerical approximations denoted by $\rho^{n+c,k},$ that approximates $\rho(t^{n}+c \Delta t)$ to $k$-th order accuracy.  

Denote by $U^{(k)}$ an approximation to the flow operator $U$ with at least $k$-th order accuracy. Then, to obtain the $(k+1)$-th order accuracy, we apply the following iterative update:
\begin{align}
    \rho^{n+1,k+1} =& U^{(k+1)}(t^{n+1},t^n) \rho^n U^{(k+1)}(t^{n+1},t^n) ^\dagger \nonumber\\
    &+\Delta t \sum_{j=1}^{N_q} \omega_j U^{(k)}(t^{n+1},t^{n+c_j\Delta t}) \mathcal{L}_L \rho^{n+c_j,k} U^{(k)}(t^{n+1},t^{n+c_j\Delta t}) ^\dagger. \label{eq:NPI_general}
\end{align}
Here  $\omega_j $ are the quadrature weights and $c_j$ define the quadrature nodes for a $(k+1)$-th approximation of the integral in \eqref{eqn: rho_n+1}. Note that it is the $\Delta t$ in front of the sum that allows us to use one order less accurate flow operators and the previous (less accurate) iterate $\rho^{n+c_j,k}$ and still obtain $\rho^{n+1,k+1}$ of order of accuracy $(k+1)$.  A standard local truncation error analysis can be easily carried out which will show $(k+1)$-th order of accuracy using an induction argument. To save space, we do not include the analysis in this paper.

 We  note that our method by NPI is naturally an embedded integrator that can be used to perform time-adaptive calculations. 
 The intuition behind NPI   is similar to SDC, which is based on a successive improvement of approximations, with the distinction that the SDC uses the same quadrature rules at different correction levels. However, the SDC cannot be directly applied here because it sums terms with both negative and positive signs and therefore loses the CP property.

\subsection{Examples of NPI Schemes Based on Gaussian Quadrature}
Given the requirement of positive quadrature weights it is natural to consider Gaussian quadrature formulae. Then a first-order approximation can be expressed as follows    
\begin{equation}
\label{eq:firstordernnpi}
\rho^{n+1,1} = U^{(1)}(t^{n+1},t^n) \rho^n U^{(1)}(t^{n+1},t^n) ^\dagger + \Delta t U^{(1)}(t^{n+1},t^n) \mathcal{L}_L \rho^n U^{(1)}(t^{n+1},t^n) ^\dagger.
\end{equation}

A second-order method is obtained by approximating the integral by the trapezoidal rule 
 \begin{align}
 \label{eq:npi2}
\rho^{n+1,2} =& U^{(2)}(t^{n+1},t^n) \rho^n U^{(2)}(t^{n+1},t^n) ^\dagger \\
&+ \frac{\Delta t}{2}\Big( U^{(1)}(t^{n+1},t^{n+1}) \mathcal{L}_L \rho^{n+1,1} U^{(1)}(t^{n+1},t^{n+1})^\dagger +  U^{(1)}(t^{n+1},t^n) \mathcal{L}_L \rho^n U^{(1)}(t^{n+1},t^n)^\dagger \Big) \nonumber \\
= & U^{(2)}(t^{n+1},t^n) \rho^n U^{(2)}(t^{n+1},t^n) ^\dagger + \frac{\Delta t}{2} \Big( \mathcal{L}_L \rho^{n+1,1}+  U^{(1)}(t^{n+1},t^n) \mathcal{L}_L \rho^n U^{(1)}(t^{n+1},t^n) ^\dagger\Big). \nonumber
\end{align}
Note here the nesting, with the first-order accurate approximate solution entering inside the quadrature approximation. The resulting method resembles the exponential trapezoidal method. 

Another option for a second order accurate method is to use midpoint integration rule, which will give
\begin{align}
    \rho^{n+1,2} = & U^{(2)}(t^{n+1},t^n) \rho^n U^{(2)}(t^{n+1},t^n) ^\dagger \\
     &+ \Delta t\Big(  U^{(1)}(t^{n+1},t^{n+1/2}) \mathcal{L}_L \rho^{n+1/2,1} U^{(1)}(t^{n+1},t^{n+1/2}) ^\dagger\Big), \nonumber
\end{align}
where $\rho^{n+1/2,1}$ is obtained by \eqref{eq:firstordernnpi} with time step size $\Delta t/2.$

Continuing, third-order accuracy can be obtained by using the third-order Gauss-Radau quadrature for the integral,  
\begin{align}
\label{eq:npi3}
\rho^{n+1,3} =& U^{(3)}(t^{n+1},t^n) \rho^n U^{(3)}(t^{n+1},t^n) ^\dagger \\
&+ \frac{\Delta t}{4}\Big(   3U^{(2)}(t^{n+1},t^{n+\frac{2}{3}}) \mathcal{L}_L \rho^{n+\frac{2}{3},2} U^{(2)}(t^{n+1},t^{n+\frac{2}{3}}) ^\dagger +U^{(2)}(t^{n+1},t^n) \mathcal{L}_L \rho^n U^{(2)}(t^{n+1},t^n) ^\dagger\Big),\nonumber
\end{align}
where $\rho^{n+\frac{2}{3},2} $ denotes the numerical approximation of $\rho$ at time $t_n + \frac{2}{3}\Delta t$ with second-order accuracy. This quantity can, for example, be obtained by using \eqref{eq:npi2} and changing the time step from $\Delta t$ to $\frac{2}{3}\Delta t.$

The process can be continued to any order but in this paper we stop at order four where we approximate the integral using Gauss-Legendre quadrature. The fourth-order approximation of $\rho(t^{n+1})$ we use is given by 
\begin{align*}
    \rho^{n+1,4} &= U^{(4)}(t^{n+1},t^n) \rho^n U^{(4)}(t^{n+1},t^n) ^\dagger\\
& + \frac{\Delta t}{2}\Big(   U^{(3)}(t^{n+1},t^{n+c_1\Delta t}) \mathcal{L}_L \rho^{n+c_1,3} U^{(3)}(t^{n+1},t^{n+c_1\Delta t}) ^\dagger+\\
& \hspace{3cm} U^{(3)}(t^{n+1},t^{n+c_2\Delta t}) \mathcal{L}_L \rho^{n+c_2,3} U^{(3)}(t^{n+1},t^{n+c_2\Delta t}) ^\dagger\Big),
\end{align*}
where $c_{1,2}= \frac{3 \mp\sqrt{3}}{6}.$ Again, note  that $\rho^{n+c_1,3}, \rho^{n+c_2,3}$ can be computed from \eqref{eq:npi3} using appropriate time step sizes.

\subsection{Examples of Approximate Flow Operators} \label{sec:flow}
We now turn to the numerical approximation of the flow operators. By definition the flow operator $U(t^{n+1},t^n)$ 
is  the solution $V(t^{n+1})$ to   
$
    \frac{dV(t)}{dt} = J(t) V(t), $
at time $t=t^{n+1}$ and with initial condition given at $t=t^n$. An approximation to the flow operator can thus be found by approximately solving this equation. In principle, any numerical method for systems of ODE can be used but here we limit ourselves to a few well established methods. 
Precisely for first-order accuracy we consider the forward and backward Euler method, and for second-order accuracy, we consider the explicit/implicit midpoint method. 
For explicit methods of third and fourth order, we use Kutta's method and the classical fourth-order Runge-Kutta method, respectively. 
For the implicit versions of the third and fourth-order accurate NPI methods, we use the fourth-order implicit scheme from \cite{puzynin2000magnus} defined as follows.
\begin{equation}
\label{eq:fourthorderimf}
\begin{cases}
     (I - \frac{1}{4}d \Delta t F_n) V^{n+1/2} = (I - \frac{1}{4}\bar{d} \Delta t F_n) V^{n}\\
    (I + \frac{1}{4}\bar{d} \Delta t F_n) V^{n+1} = (I + \frac{1}{4}d \Delta t F_n) V^{n+1/2},
\end{cases}
\end{equation}
where 
\[
F_n = iJ(t_{n+1/2}) + \frac{(\Delta t)^2}{24}i \frac{d^2J(t_{n+1/2})}{dt^2} +\frac{(\Delta t)^2}{12}i \left( \frac{ dJ(t_{n+1/2})}{dt} J(t_{n+1/2}) - J(t_{n+1/2})\frac{dJ(t_{n+1/2})}{dt}  \right),
\]
and \mbox{$d=\frac{1}{\sqrt{3}}-i$}. 

In terms of notations, we denote, the $k$th-order explicit and implicit flow operators $U^{(k),\text{ex}}$ and $U^{(k),\text{im}}$, respectively. These choices will impact    properties of the overall method and this will be discussed in the following sections. 

\section{Stability Analysis} \label{sec:stability}
In this section, we  analyze the stability region of NPI scheme \eqref{eq:NPI_general} without trace renormalization. Since the discrete numerical operator \eqref{eq:NPI_general}  is only a linear CP map (but not linear CPTP map), the contractiveness properties do not automatically carry through. Therefore, the analysis of the eigenspectrum of the discrete numerical operator is important to validate the stability properties of the methods.

\subsection{Single-qubit Test Equation}
\label{sec:teste1s}
To guarantee physical relevance, instead of applying our methods to the Dahlquist test equation as is the usual   practice for standard numerical ODE methods, we use two physics based test equations originating from  a single-qubit quantum system. In this subsection, we consider the model equation with   Hamiltonian and the Lindblad operators   given by
\[
H =  \omega \bm{0&0\\0&1}, L= \frac{1}{\sqrt{T}}\bm{0&1\\0&0},
\]
respectively.  Here $\omega$ is  the transition frequency and  $T>0$ denotes the characteristic decay time. The jump operator here models the decay process. The  term
 $J(t)$ in \eqref{eqn: master_rewrite} combining coherent evolution and dissipative effects evaluates to
\[
J = -iH -\frac{1}{2} L^\dagger L = \bm{0&0\\0&-\omega i - \frac{1}{2T}}.
\]
Denote $\alpha_J =  -\omega i - \frac{1}{2T}$, the superoperator of this equation is then 
\[
\bm{0&0&0&\frac{1}{T}\\
0 &\overline{\alpha_J}&0&0\\
0&0&\alpha_J&0\\
0&0&0& - \frac{1}{T}},
\]
with eigenvalues $0,\overline{\alpha_J},\alpha_J, -\frac{1}{T}$.  

We next derive amplification matrices for numerical approximations of different orders of accuracy. Specifically, the discussions below will illustrate how the stability of the  methods depends on the timestep $\Delta t$ and the parameters $\omega$ and $T$. 

To simplify the expressions of the amplification matrices and the associated stability conditions, we introduce:
\begin{align*}
    \alpha_{0,s} &= 1+s\Delta t\alpha_J, \quad
    \alpha_{1,s} = \frac{1}{1-s\Delta t\alpha_J}, \\
  \beta_{0,s} &=  1+s\Delta t\alpha_J +\frac{(s\Delta t)^2}{2} (\alpha_J)^2, \quad
  \beta_{1,s} =\frac{1+\frac{s\Delta t}{2}\alpha_J}{1-\frac{s\Delta t}{2}\alpha_J} ,\\
  \gamma_{0,s} &= 1+s\Delta t\alpha_J +\frac{(s\Delta t)^2}{2} (\alpha_J)^2+\frac{(s\Delta t)^3}{6} (\alpha_J)^3,\\
  \gamma_{1,s} & = 
  \frac{( 1+i\frac{s\Delta t}{4}d \alpha_J)( 1-i\frac{s\Delta t}{4}\bar{d} \alpha_J)}{(1+i\frac{s\Delta t}{4}\bar{d} \alpha_J)( 1-i\frac{s\Delta t}{4}d \alpha_J)}  ,\\
  \delta_{0,s} &= 1+s\Delta t\alpha_J +\frac{(s\Delta t)^2}{2} (\alpha_J)^2+\frac{(s\Delta t)^3}{6} (\alpha_J)^3+\frac{(s\Delta t)^4}{24}(\alpha_J)^4,\\
\delta_{1,s}&=\gamma_{1,s},
\end{align*}
where (again) $d=\frac{1}{\sqrt{3}}-i$.  

\subsubsection{First-order Schemes}
A first-order approximation using   forward Euler scheme for the flow operator $U^{(1),\text{ex}} = I + \Delta t J$,   can be formulated as:
\begin{equation}
\label{eq:firstnpiex}
\rho^{n+1,1} = (I+\Delta t J) \rho^n (I+\Delta t J)^\dagger + \Delta t (I+\Delta t J) L \rho^n L^\dagger  (I+\Delta t J)^\dagger.
\end{equation}
Similarly, using the backward Euler scheme for the flow operator, $U^{(1),\text{im}} = (I - \Delta t J)^{-1},$ we obtain:
\begin{equation}
\label{eq:firstnpiim}
\rho^{n+1,1} =(I-\Delta t J)^{-1} \rho^n ((I-\Delta t J)^{-1})^\dagger + \Delta t (I-\Delta t J)^{-1} L \rho^n L^\dagger  ((I-\Delta t J)^{-1})^\dagger.
\end{equation}
The amplification matrix for both the first order accurate methods can be expressed as
\[
\bm{ 1 &0&0&\frac{\Delta t}{T}\\ 0&\overline{\alpha_{i,1}}& 0&0\\0&0&\alpha_{i,1}&0\\ 0&0&0&\left|\alpha_{i,1}\right|^2},
\]
where then
\[
i= \begin{cases} 
0 , & \text{if } U^{(1)}=U^{(1),\text{ex}}, \\
 1,& \text{if } U^{(1)}=U^{(1),\text{im}},
\end{cases}
\]
depending what method is used. Clearly, there is a stationary solution corresponding to the eigenvalue $1.$ The stability condition for the explicit scheme \eqref{eq:firstnpiex} is
\[
( 1-\frac{\Delta t}{2T})^2+(\omega \Delta t)^2\le1.
\]
While for the implicit scheme \eqref{eq:firstnpiim}, the stability condition becomes
\[
\left| \frac{1}{1+\Delta t\omega i + \frac{\Delta t}{2T}} \right|\le1,
\]
which implies that \eqref{eq:firstnpiim} is unconditionally stable.

\subsubsection{Second-order Schemes}
The second-order approximation of the solution using the midpoint approximation for the integral is defined by
\begin{equation}
\label{eq:secondnpix}
    \rho^{n+1,2} = U^{(2)}(t^{n+1},t^n) \rho^n U^{(2)}(t^{n+1},t^n) ^\dagger  + \Delta t\Big(  U^{(1)}(t^{n+1},t^{n+\frac{1}{2}}) \mathcal{L}_L \rho^{n+\frac{1}{2},1} U^{(1)}(t^{n+1},t^{n+\frac{1}{2}}) ^\dagger\Big),
\end{equation}
    where  $\rho^{n+\frac{1}{2},1}$ is obtained by \eqref{eq:firstnpiex} or \eqref{eq:firstnpiim} with time step $\Delta t/2,$ denoted by $\rho^{n+1/2,1} = \rho^{n+1/2,1,\text{ex}}$ or $\rho^{n+1/2,1} = \rho^{n+1/2,1,\text{im}},$ respectively.
We consider explicit second-order scheme for the flow 
\[
U^{(2)}=U^{(2),\text{ex}} = I + \Delta t J + \frac{(\Delta t)^2}{2}J^2,
\]
as well as the implicit version 
\[
U^{(2)}=U^{(2),\text{im}} = (I - \frac{\Delta t}{2}J)^{-1}(I + \frac{\Delta t}{2}J).
\]
There are a total of eight combinations of implicit/explicit combination possible for the three flow operators involved. 
The amplification matrices   all take the form
\[
\bm{ 1 &0&0&\frac{\Delta t}{T} \left|\alpha_{j,\frac{1}{2}}\right|^2 \\ 0&\overline{\beta_{i,1}}& 0&0\\0&0&\beta_{i,1}&0\\ 0&0&0&\left|\beta_{i,1}\right|^2},
\]
where
\[
i= \begin{cases} 
0 , & \text{if } U^{(2)}=U^{(2),\text{ex}}, \\
 1,& \text{if } U^{(2)}=U^{(2),\text{im}},
\end{cases}
\qquad 
j = \begin{cases} 
0 , & \text{if } \rho^{n+1/2,1} = \rho^{n+1/2,1,\text{ex}}, \\
 1,& \text{if }  \rho^{n+1/2,1} = \rho^{n+1/2,1,\text{im}}.
\end{cases}
\]

Note that the stability of these numerical schemes only depend on the second-order flow operator. In particular, as long as $U^{(2)}=U^{(2),\text{im}}$ the scheme is unconditionally stable. The stability condition  when $U^{(2)}=U^{(2),\text{ex}}$ is
\[
\left |\frac{1}{2} \left(-\Delta t\omega i - \frac{\Delta t}{2T}+1\right)^2 +\frac{1}{2} \right |\le1.
\]

\subsubsection{Third- and Fourth-order Schemes}
Due to the observation made in the previous subsection, for simplicity, we will use explicit flow operators except for the highest order one. This means we fix $U^{(1)}=U^{(1),\text{ex}},U^{(2)}=U^{(2),\text{ex}}.$ For the third-order scheme, we consider  
\[
U^{(3)}=U^{(3),\text{ex}}= I + \Delta t J + \frac{(\Delta t)^2}{2}J^2+\frac{(\Delta t)^3}{6}J^3,
\] 
and
\[
U^{(3)}=U^{(4),\text{im}}=(I +  i\frac{\Delta t}{4}\overline{d} J)^{-1}(I + i\frac{\Delta t}{4} d  J)\,
(I - i\frac{\Delta t}{4}d J)^{-1}(I - i\frac{\Delta t}{4}\overline{d}  J). 
\]  
For the fourth-order scheme, we consider 
\[
U^{(4)}=U^{(4),\text{ex}}=I + \Delta t J + \frac{(\Delta t)^2}{2}J^2+\frac{(\Delta t)^3}{6}J^3+\frac{(\Delta t)^4}{24}J^4,
\]
and $U^{(4)}=U^{(4),\text{im}}.$

 For the third-order scheme, the amplification matrix is 
\[
\bm{ 1 &0&0&\frac{\Delta t}{4T} +3\frac{\Delta t}{4T}\left|\beta_{0,2/3}\right|^2 \\ 0&\overline{\gamma_{i,1}}& 0&0\\0&0&\gamma_{i,1}&0\\ 0&0&0&\left|\gamma_{i,1}\right|^2},
\]
where
\[
i = \begin{cases}
     0& \text{ if } U^{(3)} = U^{(3), \text{ex}}, \\1 
     & \text{ if }U^{(3)} = U^{(4), \text{im}} .
\end{cases}
\]
And for the fourth-order scheme, the amplification matrix is 
 \[
\bm{ 1 &0&0&\frac{\Delta t}{2T}(   \left|\gamma_{0,c_1} \right|^2+ \left|\gamma_{0,c_2} \right|^2  )\\ 0&\overline{\delta_{i,1}}& 0&0\\0&0&\delta_{i,1}&0\\ 0&0&0&\left|\delta_{i,1}\right|^2},
\]
where 
\[
i = \begin{cases}
     0& \text{ if } U^{(4)} = U^{(4), \text{ex}}, \\1 
     & \text{ if }U^{(4)} = U^{(4), \text{im}} ,
\end{cases}
\qquad \text{and }\qquad c_{1,2}= \frac{3 \mp\sqrt{3}}{6}.
\]
The stability conditions for third- and fourth-order explicit schemes are 
\[\left| 1-\Delta t \omega i - \frac{\Delta t}{2T}+\frac{1}{2}(-\Delta t \omega i - \frac{\Delta t}{2T})^2+\frac{1}{6}(-\Delta t \omega i - \frac{\Delta t}{2T})^3\right| \le 1,\] 
 and \[\left| 1-\Delta t \omega i - \frac{\Delta t}{2T}+\frac{1}{2}(-\Delta t \omega i - \frac{\Delta t}{2T})^2+\frac{1}{6}(-\Delta t \omega i - \frac{\Delta t}{2T})^3+ \frac{1}{24}(-\Delta t \omega i - \frac{\Delta t}{2T})^4\right| \le 1,\]  respectively.   
For the implicit methods, because the flow operator \eqref{eq:fourthorderimf} corresponds to  $(2,2)$-Pad\'{e} approximation in this case, the method is unconditionally stable \cite{wanner1996solving}.


\begin{figure}[htb]
\centering
\includegraphics[width=0.6\textwidth, trim=0cm 0cm 0cm 0cm, clip=true]{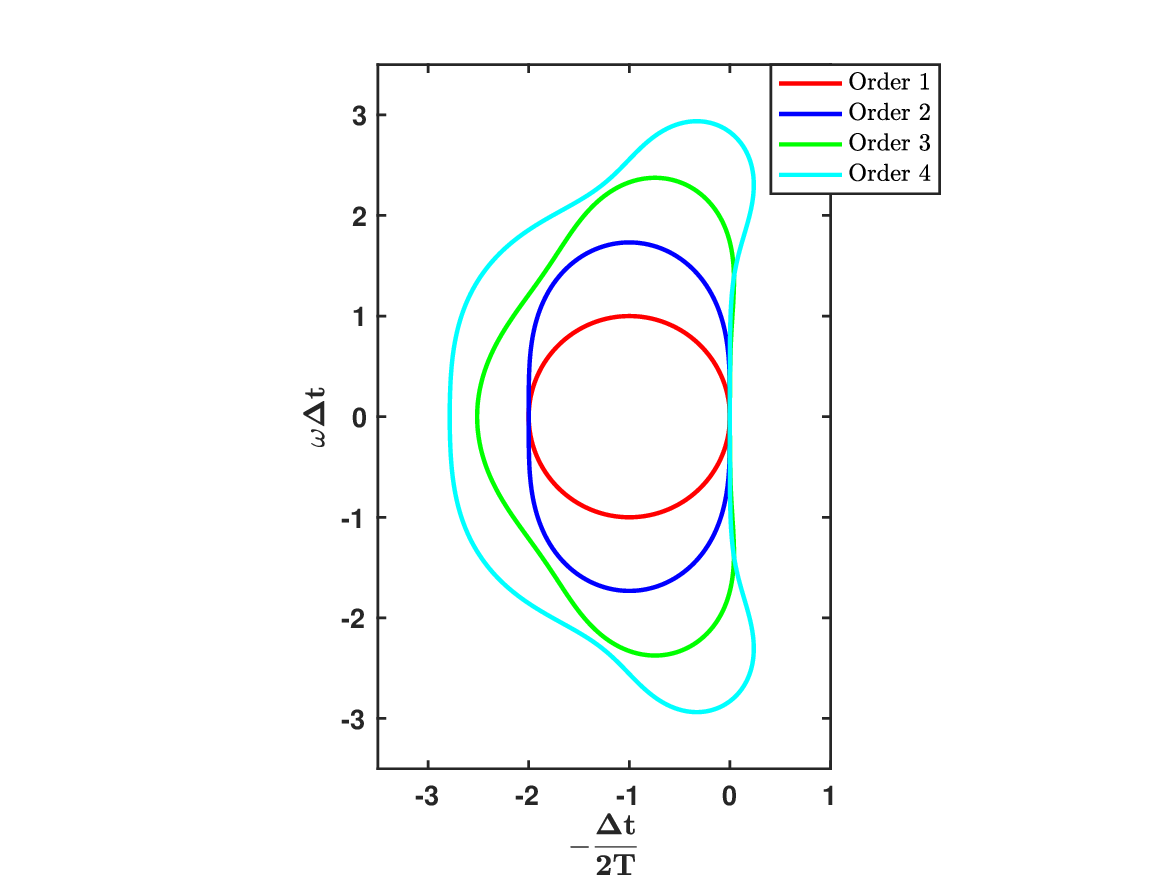}
\caption{ Stability region of the     schemes corresponding to schemes with  highest order flow operator taken as explicit schemes of order 1, 2, 3, 4 for the single-qubit test equation in Section \ref{sec:teste1s}: the interior domains bounded by the  curves and the rightmost boundary at $\Delta t=0$.}
\label{fig:stability_3_4_model_1}
\end{figure}

In conclusion, we find that for this test equation, the deciding factor for stability is the stability of the highest order flow operator. This is not surprising due to the nature of the jump operator considered. If we use implicit discretization for the highest order flow operator, the scheme will be unconditionally stable. If we use an explicit discretization for the highest order flow operator, we will have a step size restriction, which will correspond to stability region from the corresponding explicit  flow operator. For reference, we plot the stability region obtained by explicit flow operators from order 1 to order 4 in Figure \ref{fig:stability_3_4_model_1}, which is of the same shape as the standard explicit Runge-Kutta schemes.

\subsection{Single-qubit Test Equation: More General Case}
Now we consider a more general setup, where 
the Hamiltonian and the Lindblad operators of the model are given by
\[
H =  \omega \bm{0&0\\0&1}, L_1= \frac{1}{\sqrt{T_1}}\bm{0&1\\0&0},  \ \ L_2 = \frac{1}{\sqrt{T_2}}\bm{0&0\\0&1}.
\]
The jump operator $L_1$  models decay, while the $L_2$ operator models dephasing.
This will give rise to
\[
J = -iH -\frac{1}{2} L_1^\dagger L_1 -\frac{1}{2} L_2^\dagger L_2 = \bm{0&0\\0&-\omega i -\frac{1}{2T_1}-\frac{1}{2T_2}}.
\]
Denote $\alpha_J : =  -\omega i -\frac{1}{2T_1}-\frac{1}{2T_2}$. The superoperator of this equation is 
\[
\bm{0&0&0&\frac{1}{T_1}\\0 &\overline{\alpha_J}&0&0\\0&0&\alpha_J&0\\0&0&0& -\frac{1}{T_1}},
\]
with eigenvalues $0,\overline{\alpha_J},\alpha_J, -\frac{1}{T_1}$. 
We carry out the same stability analysis as above for first through fourth-order schemes. 

The amplification matrix for the first-order schemes is
\[
\bm{ 1 &0&0&\frac{\Delta t}{T_1} \\ 0&\overline{\alpha_{i,1}}& 0&0\\0&0&\alpha_{i,1}&0\\ 0&0&0&\frac{\Delta t+T_2}{T_2}\,|\alpha_{i,1}|^{2}},\]
 where \begin{align}\label{(ii)-indices}
 i=\begin{cases}
    0 & \text{ if } U^{(1)} = U^{(1),\text{ex}},  \\
    1 & \text{ if } U^{(1)} = U^{(1),\text{im}},
\end{cases}
\end{align}
which gives the stability condition 
\[
\; 
|\alpha_{i,1}| \le \sqrt{\frac{T_2}{\Delta t +T_2}}.
\]
We note that this is a more stringent condition than the one obtained in the previous subsection. In particular, the implicit scheme is no longer unconditionally stable. However, if the dephasing effects are weak, i.e. $T_2$ is large, the impact on stability condition will be small.

The amplification matrix for the second-order scheme based on the midpoint rule is
\[
\bm{ 1 &0&0&\frac{\,\Delta t\,(\Delta t+2T_2)}{2\,T_1\,T_2}|\alpha_{j,\frac{1}{2}}|^{2} \\ 0&\overline{\beta_{i,1}}& 0&0\\0&0&\beta_{i,1}&0\\ 0&0&0&\frac{\Delta t\,(\Delta t+2T_2)}{2\,T_2^{2}}\,
\bigl|\alpha_{j,\tfrac{1}{2}}\bigr|^{2}\,
\bigl|\alpha_{l,\tfrac{1}{2}}\bigr|^{2}
+\lvert \beta_{i,1}\rvert^{2}},
\]
where
\begin{align}\label{(i,j,l)-indices}
i =\begin{cases}
    0 & \text{ if } U^{(2)} = U^{(2),\text{ex}}\\
     1 & \text{ if } U^{(2)} = U^{(2),\text{im}}
\end{cases},\quad
j=\begin{cases}
    0 & \text{ if } \rho^{n+1/2,1} = \rho^{n+1/2,1,\text{ex}}\\
     1 & \text{ if } \rho^{n+1/2,1} = \rho^{n+1/2,1,\text{im}}
\end{cases}, \quad
l=\begin{cases}
    0 & \text{ if } U^{(1)} = U^{(1),\text{ex}}\\
     1 & \text{ if } U^{(1)} = U^{(1),\text{im}}
\end{cases}.   
\end{align}
The stability condition for this second-order scheme is therefore
\[
\;
\left|\frac{\Delta t\,(\Delta t+2T_2)}{2\,T_2^{2}}\,
\bigl|\alpha_{j,\tfrac{1}{2}}\bigr|^{2}\,
\bigl|\alpha_{l,\tfrac{1}{2}}\bigr|^{2}
+\lvert \beta_{i,1}\rvert^{2}\right| \le 1.
\]

Continuing, the amplification matrix for the third-order scheme is
\[
\bm{1 & 0 & 0 & \frac{2(\Delta t)^2 \,(\Delta t + 3T_{2})\,}{12\,T_{1}\,T_{2}^{2}}\,|\alpha_{0,\tfrac{1}{3}}|^{4}
\;+\; \frac{\Delta t}{4\,T_{1}\,}\bigl(1 + 3\,|\beta_{0,\tfrac{2}{3}}|^{2}\bigr)\\
0 & \overline{\gamma_{i,1}} & 0 & 0 \\
0 & 0 & \gamma_{i,1} & 0 \\
0 & 0 & 0 & Y^{(3)}},
\]
where
\begin{align*}
Y^{(3)}
&=
\frac{\Delta t}{4\,T_{2}}\,|\beta_{0,1}|^{2}
\;+\;
\frac{(\Delta t)^{2}(\Delta t+3T_{2})}{6\,T_{2}^{3}}\,|\alpha_{0,\tfrac{1}{3}}|^{4}\,|\beta_{0,\tfrac{1}{3}}|^{2}
\;+\;
\frac{3\,\Delta t}{4\,T_{2}}\,|\beta_{0,\tfrac{1}{3}}|^{2}\,|\beta_{0,\tfrac{2}{3}}|^{2}
\;+\;
|\gamma_{i,1}|^{2},\notag
\\
   i&= \begin{cases}
    0 & \text{ if } U^{(3)} = U^{(3),\text{ex}},\\
     1 & \text{ if } U^{(3)} = U^{(4),\text{im}}.
\end{cases}
\end{align*}
The stability condition for the third-order scheme is
\[
\left| \frac{\Delta t}{4\,T_{2}}\,|\beta_{0,1}|^{2}
\;+\;
\frac{(\Delta t)^{2}(\Delta t+3T_{2})}{6\,T_{2}^{3}}\,|\alpha_{0,\tfrac{1}{3}}|^{4}\,|\beta_{0,\tfrac{1}{3}}|^{2}
\;+\;
\frac{3\,\Delta t}{4\,T_{2}}\,|\beta_{0,\tfrac{1}{3}}|^{2}\,|\beta_{0,\tfrac{2}{3}}|^{2}
\;+\;
|\gamma_{i,1}|^{2}\right|\le 1.
\]
The amplification matrix for the fourth-order scheme is 
\[
\bm{1&0&0&\frac{\Delta t}{24\,T_{1}\,T_{2}^{3}}X\\
0&\overline{\delta_{i,1}}&0&0\\
0&0&\delta_{i,1}&0\\0&0&0&
 |\delta_{i,1}|^{2}
+ \frac{\Delta t}{24\,T_{2}^{2}}Y^{(4)}},
\]
where 
\begin{align}
    X&=\sum_{k=1}^2\Big(
3\,\Delta t\,T_2^{2}\,c_k\,\lvert \beta_{0,c_k}\rvert^{2}
+2\,(\Delta t)^{2}\,
   c_k^{2}(\Delta t c_k + 3T_2 ) \,\lvert \alpha_{0,\frac{1}{3}c_k}\rvert^{4}\,
   \lvert \beta_{0,\frac{1}{3}c_k}\rvert^{2}\nonumber\\
   &\hspace{15mm}+9\,\Delta t\,T_2^{2}\,c_k\,
   \lvert \beta_{0,\frac{1}{3}c_k}\rvert^{2}\,
   \lvert \beta_{0,\frac{2}{3}c_k}\rvert^{2}
+12\,T_2^{3}\,\lvert \gamma_{0,c_k}\rvert^{2}
\Big),  \nonumber
\\
Y^{(4)} =& 
\sum_{k=1}^2\Big(
3\,\Delta t\, c_k\,\lvert \beta_{0,c_k}\rvert^{2}\,\lvert \gamma_{0,1-c_k}\rvert^{2}
+\frac{2(\Delta t)^{2}}{T_2^{2}} \,
  c_k^{2}\,
  (c_k\,\Delta t+3\,T_2)\,
  \lvert \alpha_{0,\frac{1}{3}c_k}\rvert^{4}\,
  \lvert \beta_{0,\frac{1}{3}c_k}\rvert^{2}\,
  \lvert \gamma_{0,1-c_k}\rvert^{2} \nonumber
\\
&\hspace{13mm}+
9\,\Delta t\,c_k\,
  \lvert \beta_{0,\frac{1}{3}c_k}\rvert^{2}\,
  \lvert \beta_{0,\frac{2}{3}c_k}\rvert^{2}\,
  \lvert \gamma_{0,1-c_k}\rvert^{2} 
+12\,T_2\,
  \lvert \gamma_{0,c_k}\rvert^{2}\,
  \lvert \gamma_{0,1-c_k}\rvert^{2}
\Big), \nonumber
 \\
 i&= \begin{cases}
    0 & \text{ if } U^{(4)} = U^{(4),\text{ex}},\\
     1 & \text{ if } U^{(4)} = U^{(4),\text{im}}.
\end{cases} 
\label{i-index in U4}
\end{align}
The stability condition for the fourth-order scheme is therefore
\[
\left|\delta_{i,1}^{2} + \tfrac{\Delta t}{24 T_{2}^{2}}\,Y^{(4)}\right| \le 1.
\]
As we can see, in all cases, the stability condition is more stringent than the ones derived in the  previous subsection.

\section{Low-rank Truncation}\label{sec:lowrank}
In our prior work \cite{appelo2024kraus}, we performed low-rank truncation by the truncated SVD, and proved that this operation is a Kraus operator. Here we use the same approach. Let $r(t)$ be the rank of the $d \times d$ density matrix and $V(t)$ be its factor, i.e. $\rho(t) = V(t)V(t)^\dagger$.  The low-rank method will store and evolve $V(t),$ which is efficient  when  $r(t) \ll d.$ 

In particular,  the numerical solution $\rho^{n+1,k+1}$ from the $k+1$-th order scheme, as described in \eqref{eq:NPI_general}, can be rewritten as 
\[
\rho^{n+1,k+1} = V^{n+1,k+1} (V^{n+1,k+1})^\dagger,
\]
where 
\begin{align}
\label{eq:trunckp1}
    V^{n+1,k+1} &= \left[U^{(k+1)}(t^{n+1},t^n)V^{n,k+1} \hspace{1cm} \cdots~~~ \sqrt{\omega_j\Delta t}  U^{(k+1)}(t^{n+1},t^{n+c_j})L_{\alpha} V^{n+c_j,k} ~~~\cdots\right].
\end{align}
Here $V^{n+1,k+1}$ needs to be truncated to maintain low-rank. This can be done by the truncated SVD with a specified tolerance. As suggested in \cite{appelo2024kraus}, a natural choice for the tolerance is to take it to be proportionate to the local truncation error. Precisely we use the truncation tolerance is $ (\kappa \Delta t)^{k+2}$ when truncating \eqref{eq:trunckp1}.  Here we chose the constant $\kappa$ to be fixed for a given computation but note that it could also be updated adaptively throughout the computation and change with each iteration. 

We note that when there are many columns in \eqref{eq:trunckp1}, a truncated SVD may not be computationally efficient. In this case, we can group columns in \eqref{eq:trunckp1} and perform truncated SVD on each group hierarchically. Due to the fact that each truncation is a CP operator, the final scheme will also be CP. Finally, after all computations have been done advancing the time to $t^{n+1},$ we perform a trace renormalization step $\rho \leftarrow \frac{\rho}{Tr(\rho)}$.

 \section{Numerical Experiments}\label{sec: numeric}
Our numerical experiments will focus on problems governed by a model describing a composite quantum system modeling multiple transmons and their electromagnetic environment \cite{ doi:10.1126/science.1231930,PhysRevLett.101.080502,  PhysRevA.76.042319}.  Here we adopt the notation from \cite{gunther2021quantum}. For transmons the Hamiltonian  consists of a time independent system Hamiltonian $H_d$ and a time dependent control Hamiltonian $H_c(t)$, that is, $H(t) = H_d + H_c(t)$. The model we use for the system Hamiltonian $H_d$ is given by
\[
H_d =\displaystyle \sum_{k=0}^{Q-1}\left( \omega_k a_k^\dagger a_k -\frac{\xi _k}{2} a_k^\dagger a_k^\dagger a_k a_k  
        +\sum_{l>k} \Big( J_{kl} (a_k^\dagger a_l+ a_k a_l^\dagger) - \xi_{kl} a_k^\dagger a_k a_l^\dagger  a_l \Big)
        \right), 
\]
where $\omega_k\geq 0$ denotes the $0\rightarrow 1 $ transition frequency and $\xi_k\geq 0$ denotes the self-Kerr coefficient of the $k$-th subsystem and $J_{kl}\geq 0$ and $\xi_{kl}\geq 0$ are the cross resonance coefficients between the  $k$-th and $l$-th subsystems. Here $a_k$ is a lowering operator acting on subsystem $k$.

The control Hamiltonian modeling the action of the external control on the quantum system is given by 
\[
H_c(t) = \sum_{k=0}^{Q-1} f^k(t) (a_k + a_k^\dagger),
\]
 where $f^k(t)$ is a real-valued time-dependent control functions. 

The Lindbladian operator $\mathcal{L}(\rho(t))$ in the model takes  the following form
 \[
 \mathcal{L}(\rho(t)) = \sum_{k=0}^{Q-1} \sum_{l=1}^2 \mathcal{L}_{lk} \rho(t) \mathcal{L}^{\dagger}_{lk}  - \frac{1}{2} \left( \mathcal{L}^{\dagger}_{lk}\mathcal{L}_{lk} \rho(t) + \rho(t) \mathcal{L}^{\dagger}_{lk}\mathcal{L}_{lk}\right),
 \]
where  $\mathcal{L}_{1 k}= \frac{1}{\sqrt{T_1^k}}a_k$ is the collapse operator modeling decay processes, and $\mathcal{L}_{2 k}= \frac{1}{\sqrt{T_2^k}}a^\dagger _k a_k$ models dephasing processes in subsystem $k$. The constants $T_1^k>0$ and $T_2^k>0$ are the half-life of the decay and dephasing processes in the $k$-th subsystem, respectively. 
 
Simulating the Lindblad equation in the lab frame can be numerically challenging when the Hamiltonian $ H(t) $ contains high-frequency terms. These terms introduce rapid oscillations, resulting in stiff differential equations that require impractically small time steps for stable integration. Therefore, as is customary, we employ the {rotating frame wave approximation}. This effectively removes many of the fast dynamics, thereby improving simulation efficiency.

The rotating frame approximation starts with the choice of the rotation frequencies $\omega_k^r$ for each oscillator. The Hamiltonians are then transformed to 
\begin{align*}
    \Tilde{H}_d(t)=&  \displaystyle \sum_{k=0}^{Q-1}\left( (\omega_k- \omega_k^r) a_k^\dagger a_k -\frac{\xi _k}{2} a_k^\dagger a_k^\dagger a_k a_k  
        -\sum_{l>k}    \xi_{kl} a_k^\dagger a_k a_l^\dagger  a_l 
        \right)\\
        &+\displaystyle \sum_{k=0}^{Q-1}  \sum_{l>k} J_{kl}\left(\cos{(\eta_{kl}t)}  (a_k^\dagger a_l+ a_k a_l^\dagger)
        +
       i \sin{(\eta_{kl}t)}  (a_k^\dagger a_l- a_k a_l^\dagger)\right),\\
       \Tilde{H}_c(t) =&  \sum_{k=0}^{Q-1} \left(p^k(\vec{\alpha}^k ,t) (a_k + a_k^\dagger) + iq^k(\vec{\alpha}^k ,t) (a_k - a_k^\dagger)\right).
\end{align*}
Here $\eta_{kl}:= \omega_k^r- \omega_l^r$ are the differences in rotational frequencies between subsystems. The real-valued laboratory frame control function $f^k(t)$ can be approximated as the following
\[
f^k(t)\approx 2p^k(t) \cos{( w_k^r)}-2q^k(t) \sin{( w_k^r)}.
\]

\begin{figure}[htb]
 \centering
\includegraphics[width=0.49\linewidth,trim={0.0cm 0.0cm 0.0cm 0.0cm},clip]{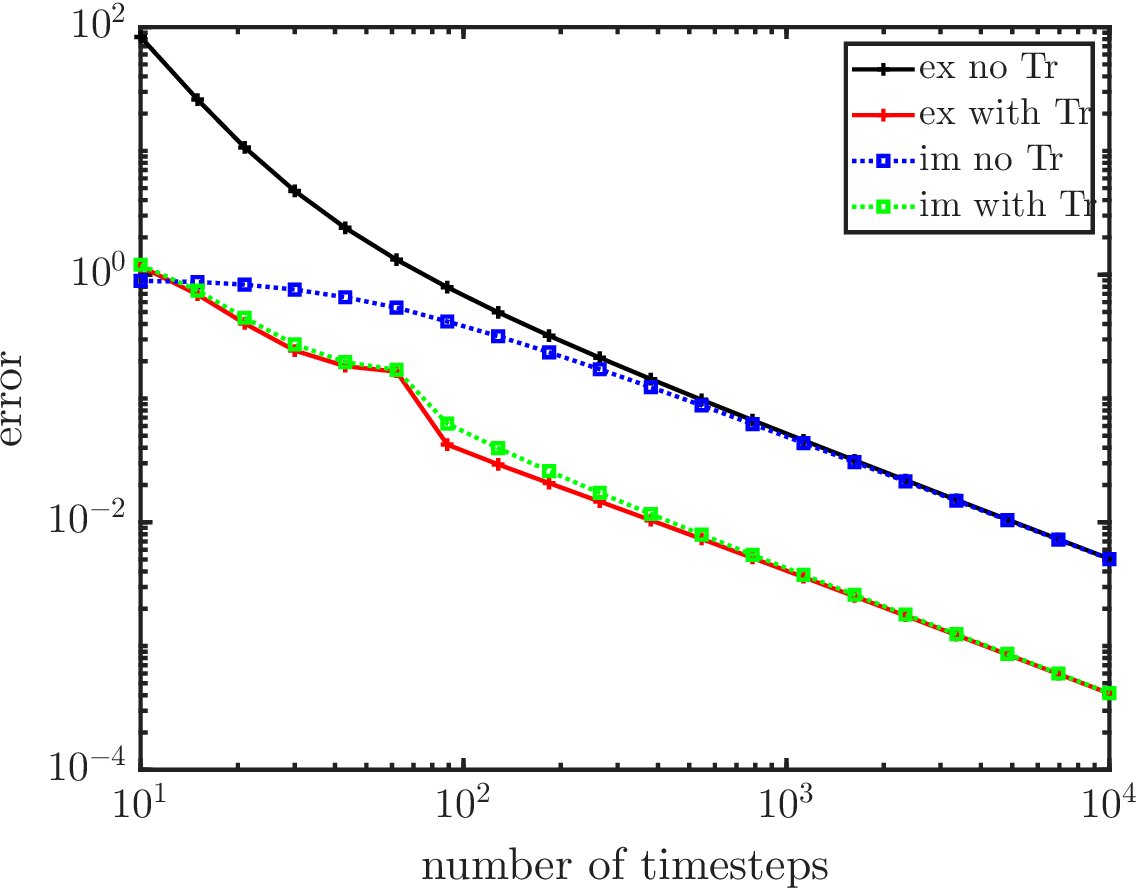}
\includegraphics[width=0.49\linewidth,trim={0.0cm 0.0cm 0.0cm 0.0cm},clip]{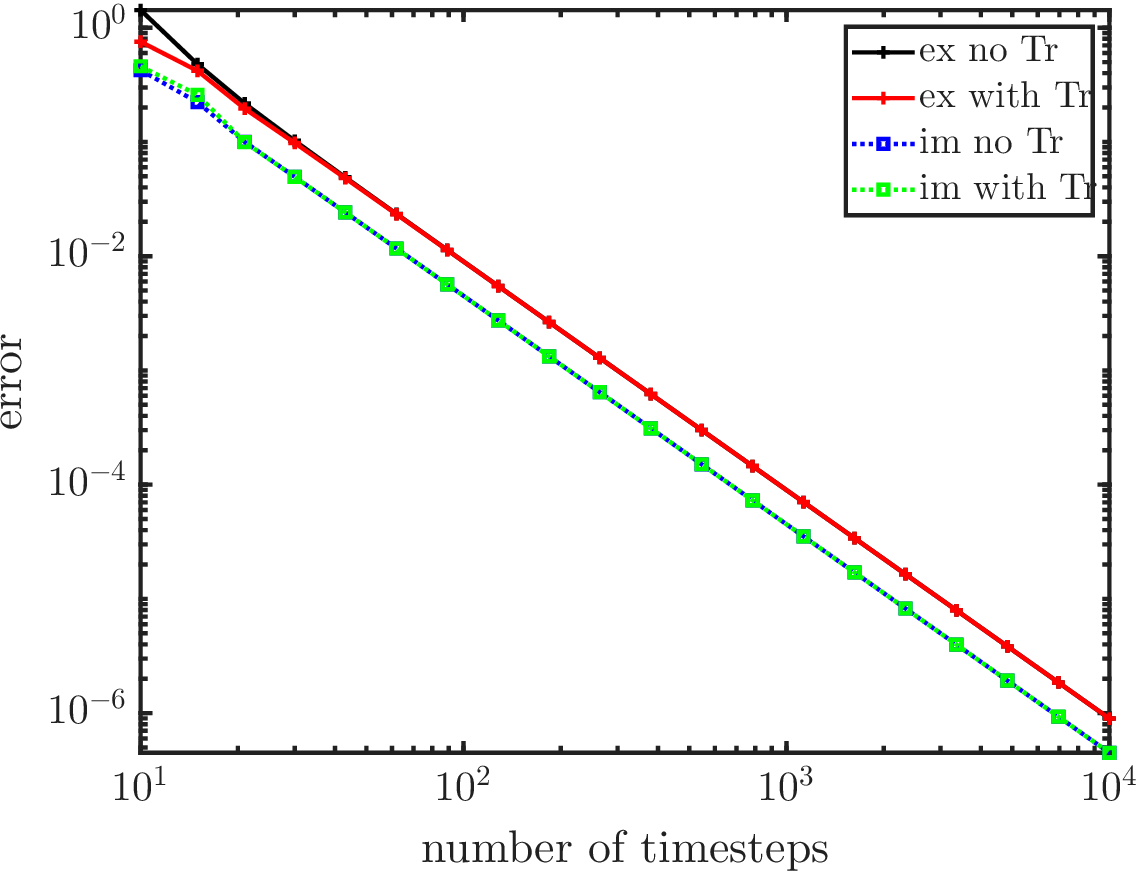}
\includegraphics[width=0.49\linewidth,trim={0.0cm 0.0cm 0.0cm 0.0cm},clip]{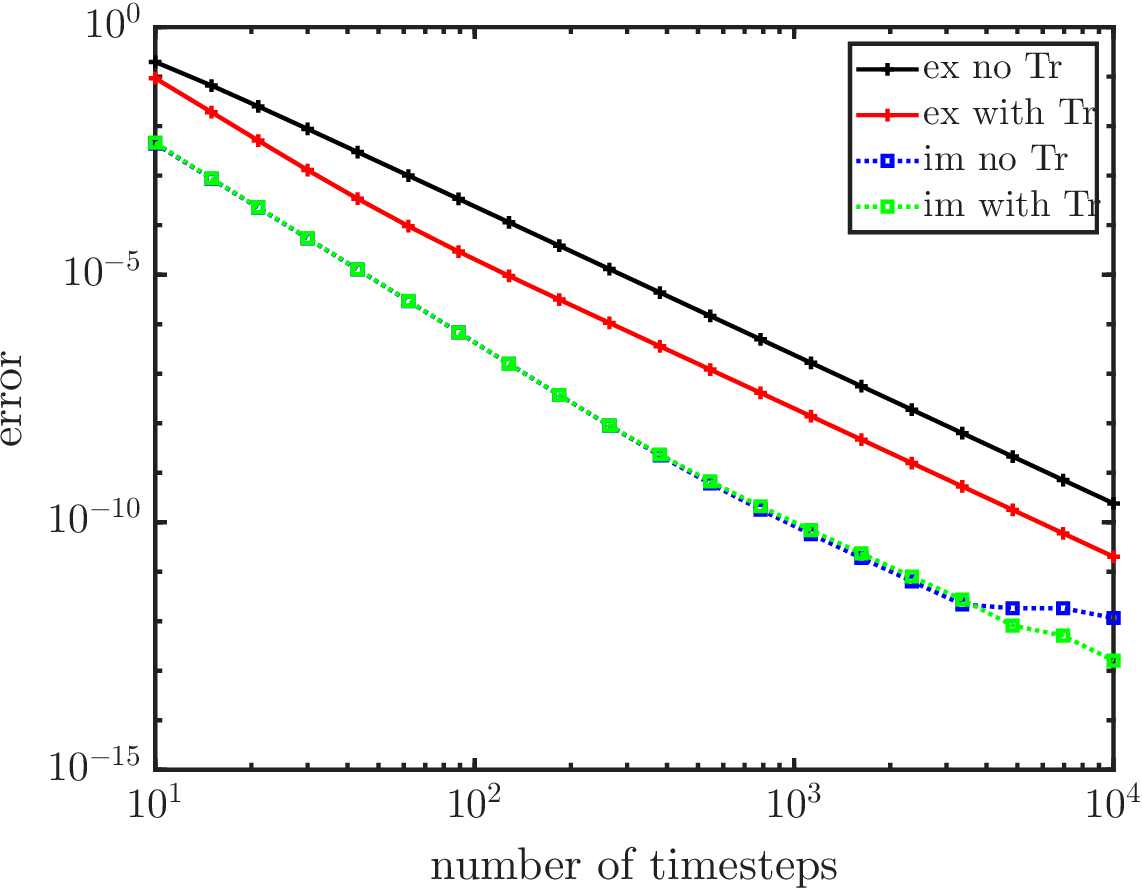}
\includegraphics[width=0.49\linewidth,trim={0.0cm 0.0cm 0.0cm 0.0cm},clip]{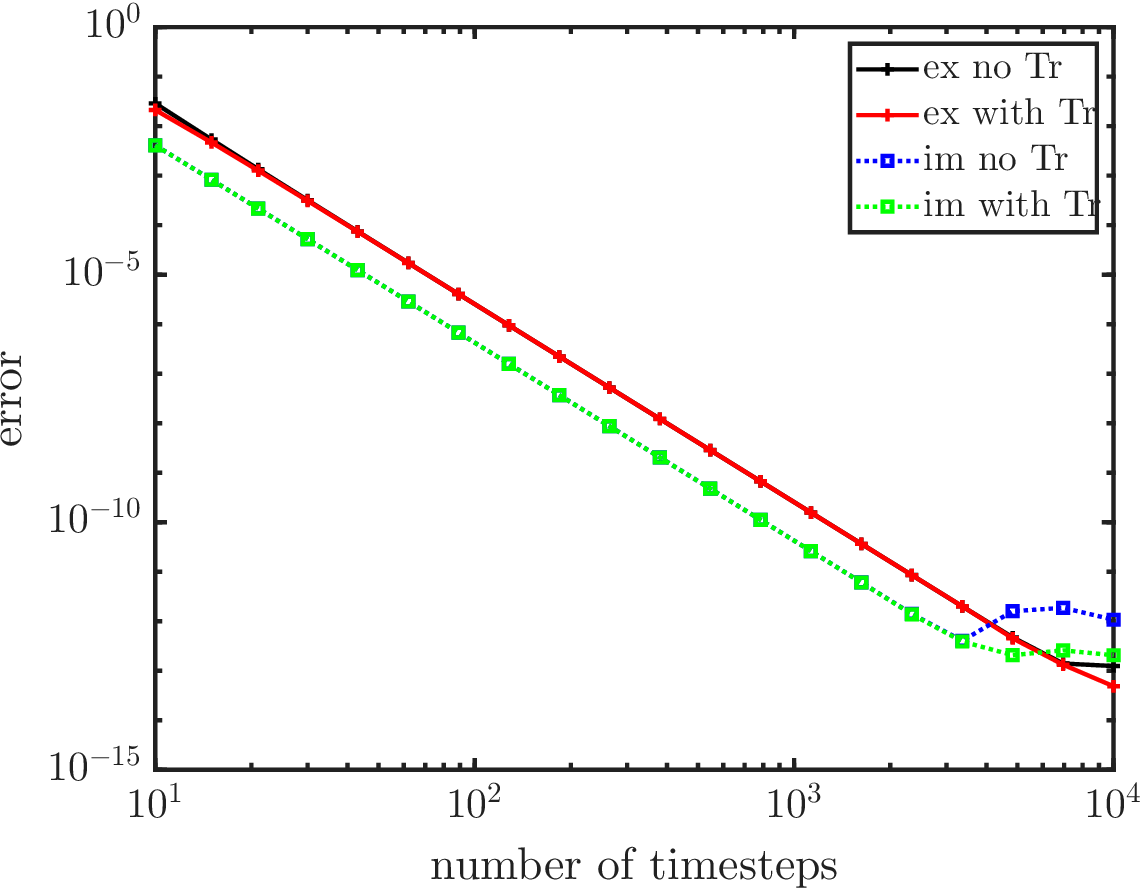}
 \caption{The graphs display the errors at the final time for methods (from top left to bottom right) of order one to four for example in Section \ref{sec:test1}. See the text for further explanation.\label{fig:2Qerr}}
 \end{figure}
\subsection{A Two-qubit System with Known Solution}
\label{sec:test1}
In this first experiment, we study the accuracy of different methods for a composite quantum system consisting of  two identical qubits with equal transition frequencies ($\omega_0 = \omega_1$). In the rotating frame, this eliminates the system Hamiltonian terms $\sum_k \omega_k a_k^\dagger a_k$. The dipole-dipole interaction is characterized by a single constant cross-resonance coefficient $J_{01} = J= 0.2$ (in non-dimensionalized units). The system Hamiltonian thus is $H = J (a_0^\dagger a_1 + a_0 a_1^\dagger)$. Here we consider Lindblad operators corresponding to decay processes with a uniform rate $\gamma = 1/T_1^0 = 1/T_1^1 = 1/50$ and exclude dephasing. Given the initial state $\rho(0) = |10\rangle\langle10|$, the exact solution at time $t$ is:
\begin{equation*}
\rho(t) = \begin{pmatrix}
1 - e^{-\gamma t} & 0 & 0 & 0 \\
0 & \frac{e^{-\gamma t}}{2}(1 - \cos (2Jt)) & -\frac{i e^{-\gamma t}}{2} \sin (2Jt) & 0 \\
0 & \frac{i e^{-\gamma t}}{2} \sin (2Jt) & \frac{e^{-\gamma t}}{2}(1 + \cos (2Jt)) & 0 \\
0 & 0 & 0 & 0
\end{pmatrix}.
\end{equation*}

\begin{table}[htb] 
 \begin{center} 
 \begin{tabular}{| l | c | c | c | c | l | c | c | c | c |} 
 \hline 
\multicolumn{5}{|c|}{first-order explicit} & \multicolumn{5}{|c|}{first-order implicit}  \\ \hline
$N_t$ & 1600 & 3200 & 6400 & 12800 & $N_t$ & 1600 & 3200 & 6400 & 12800 \\ \hline 
error & 2.6(-3) & 1.3(-3) & 6.5(-4) & 3.2(-4) & error & 2.6(-3) & 1.3(-3) & 6.5(-4) & 3.3(-4)\\ \hline 
rate &  & 0.99 & 1.00 & 1.00 & rate &  & 1.01 & 1.01 & 1.00 \\
\hline 
\hline
\multicolumn{5}{|c|}{second-order explicit} & \multicolumn{5}{|c|}{second-order implicit}  \\ \hline
$N_t$ & 200 & 400 & 800 & 1600 & $N_t$ & 200 & 400 & 800 & 1600 \\ \hline 
error & 2.2(-3) & 5.6(-4) & 1.4(-4) & 3.5(-5) & error & 1.1(-3) & 2.8(-4) & 7.0(-5) & 1.6(-5) \\ \hline 
rate &  & 2.00 & 2.00 & 2.00 & rate &  & 2.00 & 2.00 & 2.00 \\ \hline 
\hline
\multicolumn{5}{|c|}{third-order explicit} & \multicolumn{5}{|c|}{third-order implicit}  \\ \hline
$N_t$ & 45 & 90 & 180 & 360 & $N_t$ & 45 & 90 & 180 & 360 \\ \hline 
error & 2.9(-4) & 2.8(-5) & 3.4(-6) & 4.2(-7) & error & 1.1(-5) & 6.6(-7) & 4.1(-8) & 2.8(-9) \\ \hline 
rate &  & 3.36 & 3.08 & 3.00 & rate &  & 4.03 & 4.02 & 3.83 \\ \hline 
\hline
\multicolumn{5}{|c|}{fourth-order explicit} & \multicolumn{5}{|c|}{fourth-order implicit}  \\ \hline
$N_t$ & 32 & 64 & 128 & 256 & $N_t$ & 32 & 64 & 128 & 256 \\ \hline  
error & 2.4(-4) & 1.5(-5) & 9.5(-7) & 5.9(-8) & error & 4.1(-5) & 2.6(-6) & 1.6(-7) & 1.0(-8) \\ \hline  
rate &  & 3.98 & 4.00 & 4.00 & rate &  & 4.00 & 4.00 & 4.00 \\ \hline 
 \end{tabular} 
 \caption{Errors and rates of convergence at the final time for the two-qubit example with a known analytic solution. Here the trance renormalization is used. In this table $1.2(-3) = 1.2 \cdot 10^{-3}$.
 \label{tab:2Q1}}
 \end{center} 
 \end{table} 
In Figure \ref{fig:2Qerr} we display the error of the density matrix at the final time 6 in the Frobenius norm. The results are obtained using the truncation tolerance described in Section \ref{sec:lowrank} with $\kappa = 1/2$. 
 Here when we say that the method is implicit (or explicit) we mean that the approximations of the flow operators described in Section \ref{sec:flow} for each of the nested iterations are implicit (or explicit). 

Starting with the first-order methods displayed in the top left of Figure \ref{fig:2Qerr} we note that the explicit method with no trace renormalization can become unstable when using a large timestep. We further note that both the explicit methods are more accurate when the trace renormalization is applied. For the implicit methods of order two, three and four the trace renormalization does not appear to improve (or worsen) the errors. This is also the case for the second and fourth-order explicit methods. For the third-order method the trace renormalization improves the error. We don't have an explanation why the trace normalization sometimes improves the error but we note that it never deteriorates the performance of the method. In the bottom left of Figure \ref{fig:2Qerr} we can see that the implicit third-order method appears to, at first, converge at a rate close to four and then settle down to a rate of convergence that is closer to three. This is most likely due to the fact that we use the fourth-order implicit method from \cite{puzynin2000magnus} in combination with a third-order accurate quadrature.    
In Table \ref{tab:2Q1} we report errors and estimated rates of convergence for the various methods. The results are computed with the trace renormalization. The rates of convergence are consistent with the expected rates of convergence. 

\subsection{A Qudit-cavity System with Control} 

In this experiment, we consider a qutrit (subsystem 0) with $n_0=3$ energy levels coupled to a cavity (subsystem 1) with $n_1=20$ energy levels. The model parameters (adopted from \cite{gunther2021quantum}) are provided in Table \ref{table:parameters_exp_1}. 
\begin{table}[htb]
\centering 
\caption{Parameters for the qudit-cavity problem.}
\label{table:parameters_exp_1}
\begin{tabular}{|c|c|c|c|c|}
\hline  
$\omega_0/2\pi$ \SI{}{\giga\hertz} & $\omega_1/2\pi$ \SI{}{\giga\hertz} & $\xi_0/2\pi$ \SI{}{\giga\hertz} & $\xi_1/2\pi$ \SI{}{\mega\hertz} & $J_{01}$ \SI{}{\giga\hertz} \\
\hline  
$4.41666$ & $6.84081$ &  $0.23056$& $0$ &0 \\
\hline  
\hline  
$\xi_{01}$ \SI{}{\mega\hertz} & $T_{00}$ \SI{}{\micro\second} &$T_{01}$ \SI{}{\micro\second} & $T_{10}$ \SI{}{\micro\second} &$T_{11} $ \SI{}{\micro\second}\\
\hline  
 $1.176$&80&0.3892&26& $\infty$  \\
\hline  
\end{tabular}
\end{table}

To demonstrate the effect of the coupling between the two systems, we apply controls for both the qutrit and resonator on the form  
\[
p(t) = \frac{A}{2} \left(1 + \tanh (\delta ( t - \tau)) \right), \ \ q(t) = 0.
\]
We choose $A, \delta$, and $\tau$ to be \SI{10}{\mega\hertz}, \SI{-0.05}{\giga\hertz} and \SI{2000}{\nano\second} for the qutrit and \SI{15}{\mega\hertz}, \SI{-0.1}{\giga\hertz} and \SI{200}{\nano\second} for the resonator. We then simulate the problem for \SI{2.5}{\micro\second} and compute the error against a reference solution obtained using a large number of timesteps. 
\begin{figure}[htb]  
\centering
\includegraphics[width=0.60\textwidth, trim=0.0cm 0.0cm 0.0cm 0.0cm, clip=true]{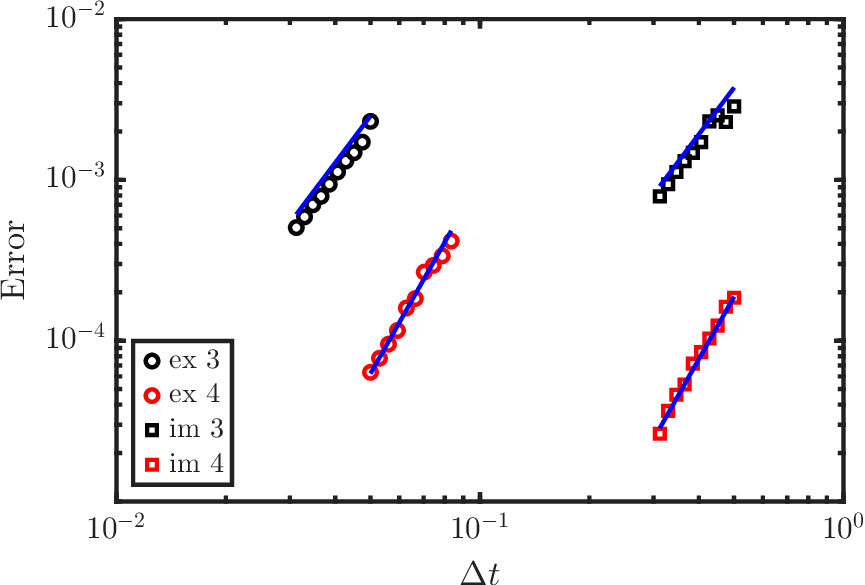}
\caption{Errors for different $\Delta t$ for the third and fourth-order implicit and explicit methods when used to compute the solution to the qubit-cavity problem.   \label{fig:QR1}}
\end{figure}

Here we only display results for the third and fourth-order methods.  We use $\kappa = 0.1$ for the implicit methods and $\kappa = 0.5$ for the explicit methods. The errors in the density matrix are measured in the Frobenius norm at the final time and displayed as a function of the number of timesteps in Figure \ref{fig:QR1}.  The smallest number of timesteps used for the explicit methods coincide with the (empirically observed) stability limit for the methods. The solid lines in the figure indicates slopes for third- and fourth-order accuracy. As can be  seen the methods perform at designed order also for problems with a  time dependent Hamiltonian. It is also clear that for this example the implicit timestepping is much more efficient in terms of the number of timesteps needed to reach a certain error. Note that the error for the explicit method is likely limited by the truncation rather than the error from the flow approximation or the error due to the quadrature. 

\begin{figure}[htb]  
\centering
\includegraphics[width=0.49\textwidth, trim=0.0cm 0.0cm 0.0cm 0.0cm, clip=true]{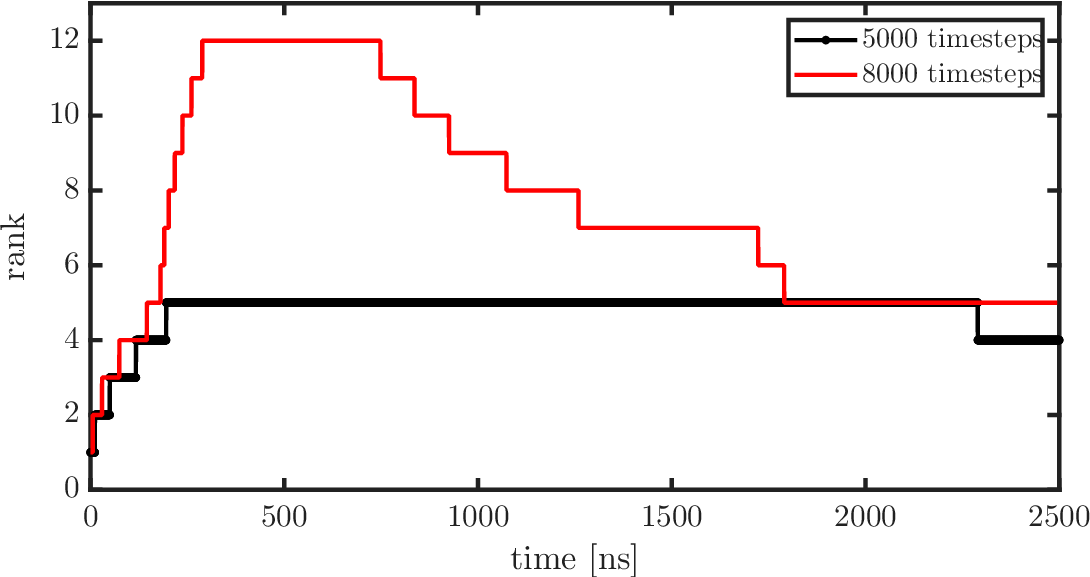}
\includegraphics[width=0.49\textwidth, trim=0.0cm 0.0cm 0.0cm 0.0cm, clip=true]{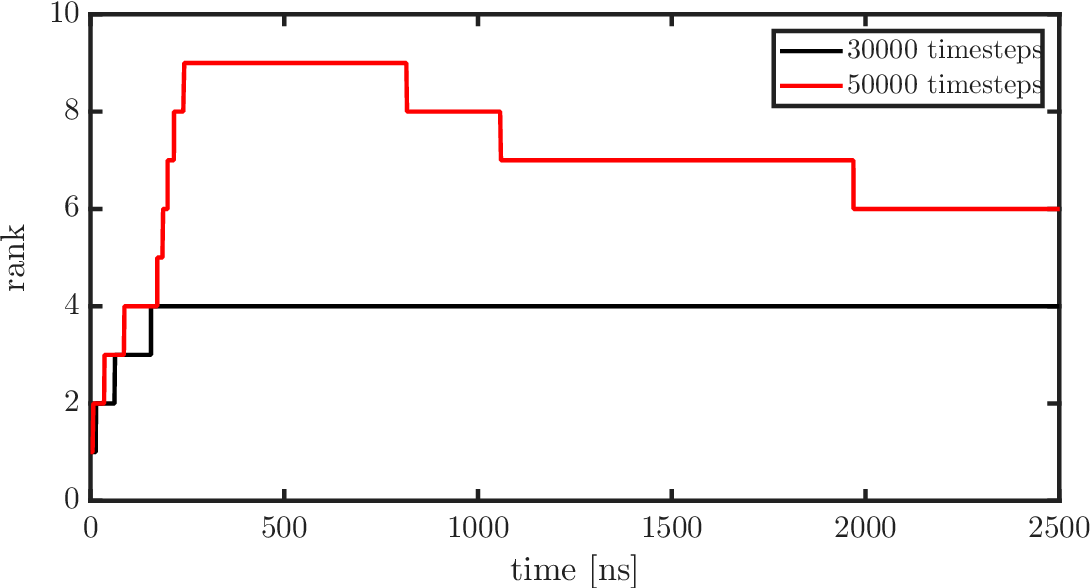}
\caption{The rank as a function of time for the qubit-cavity problem when simulated using the fourth-order implicit (left) and explicit (right) method. The number of timesteps corresponds to the largest and smallest $\Delta t$ used in the experiments producing Figure \ref{fig:QR1}. \label{fig:QR3}}
\end{figure}

In Figure \ref{fig:QR3} we display the rank as a function of time for the coarsest and finest time discretization used for the implicit and explicit fourth-order methods. As can be seen the choices $\kappa = 0.1$ and $\kappa = 0.5$ results in roughly the same level of truncation and the rank is approximately the same for the two computations. Reducing the truncation level might have reduced the error for this computation but would also have increased the cost of the method.

\begin{figure}[hbt]  
\centering
\includegraphics[width=0.49\textwidth, trim=1.8cm 0.5cm 1.5cm 0.5cm, clip=true]{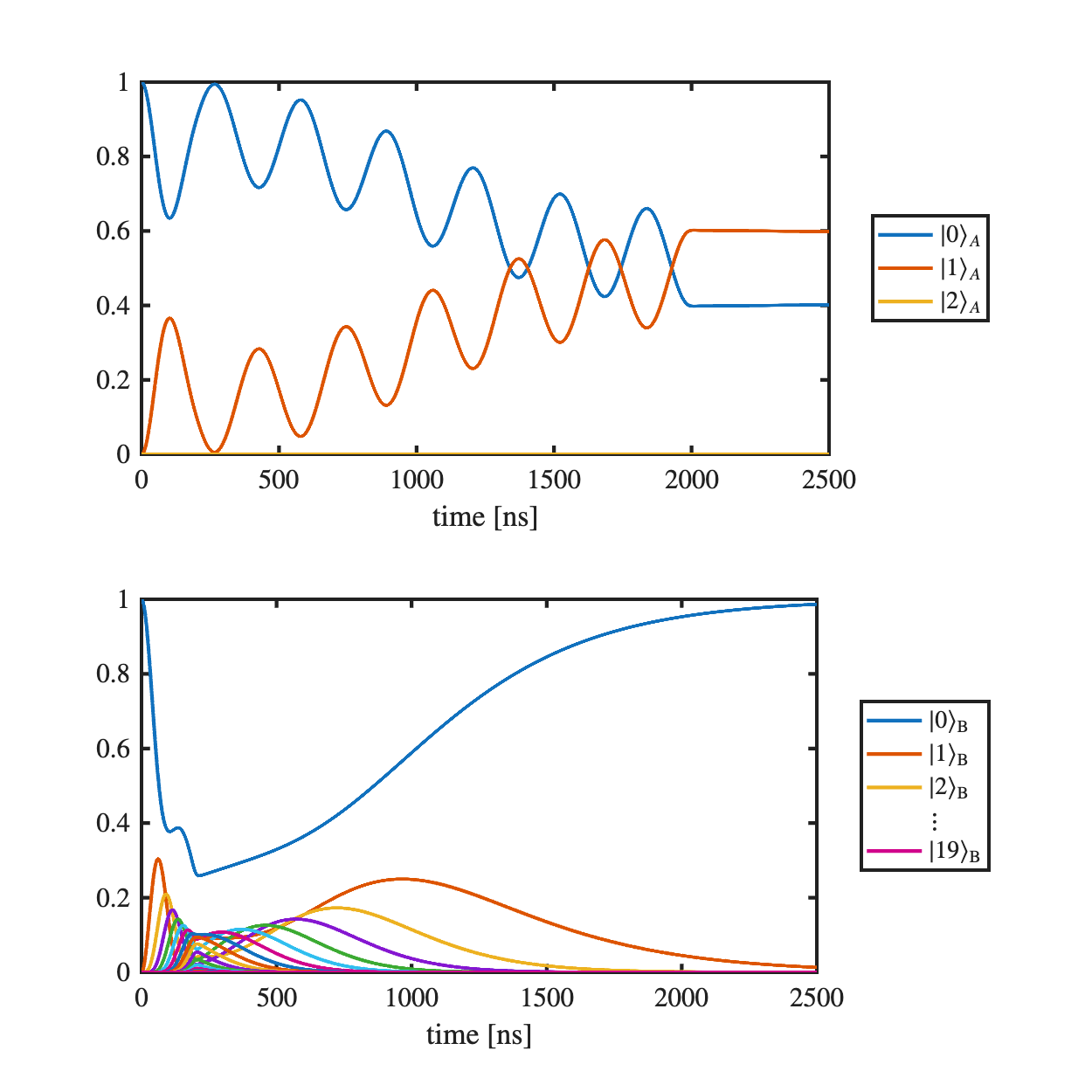}
\includegraphics[width=0.49\textwidth, trim=1.8cm 0.5cm 1.5cm 0.5cm, clip=true]{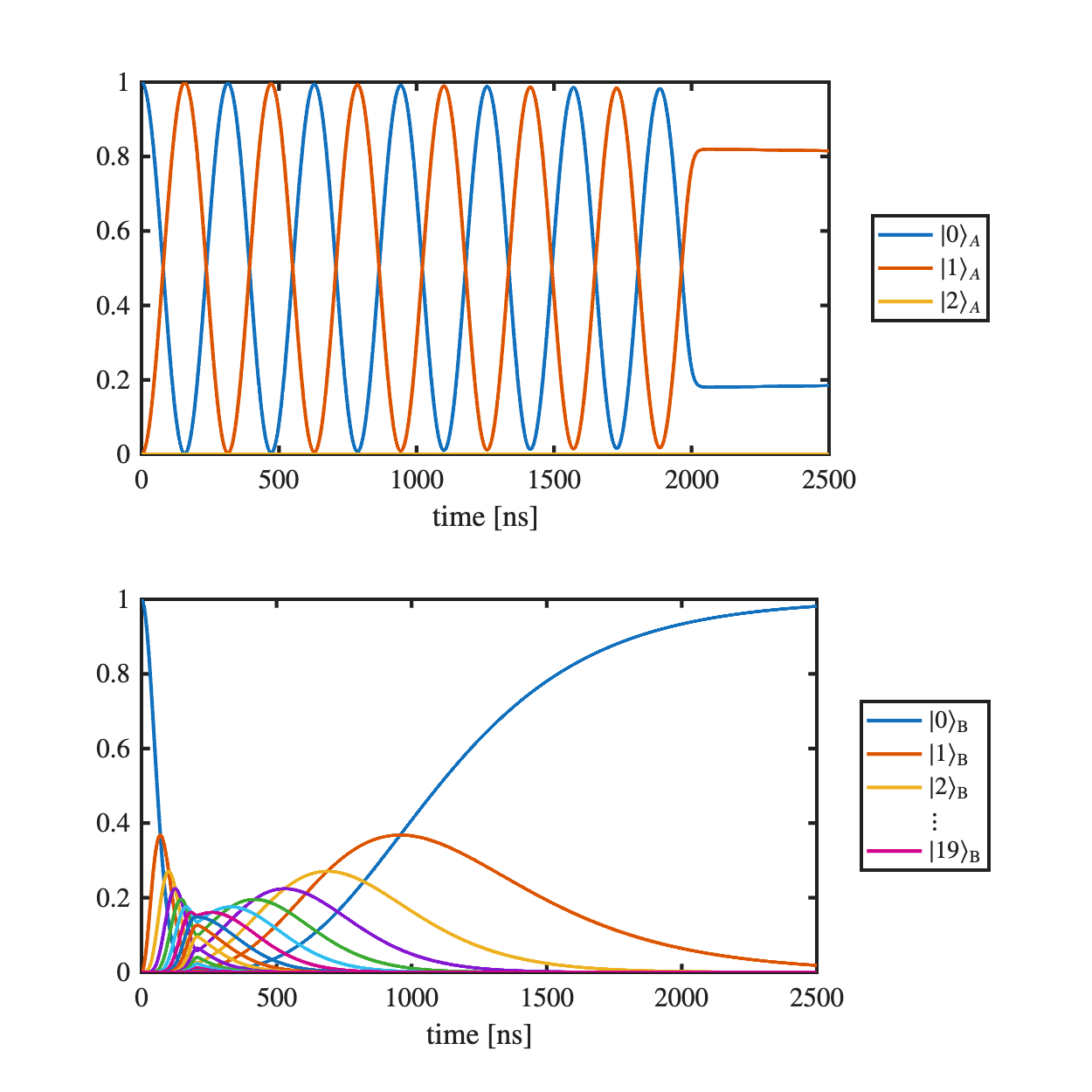}
\caption{Populations in the qubit (top, annotated by $A$) and cavity (bottom, annotated by $B$) as a function of time with coupling $\xi_{01} \neq0 $ (left) and without coupling $\xi_{01}=0$ (right).\label{fig:QR2}}
\end{figure}

Finally in Figure \ref{fig:QR2} we display the populations in the qudit and the resonator for the case with parameters as in Table \ref{table:parameters_exp_1} and as a point of reference the populations obtained with the same method and identical setup but with the coupling turned off (we set $\xi_{01} = 0$).

\subsection{Suppression of Oscillation Revival in a Jaynes-Cumming Model}
In this experiment, we consider a Jaynes-Cumming model. This model, introduced by Edwin Jaynes and Fred Cummings in 1963, explores the interaction between a two-level atom (or qubit) and a single quantized mode of an electromagnetic field, typically within a cavity environment \cite{jaynes1963comparison}. Here the model consists of a two-level atom interacting with a quantized field mode consisting of $m = 150$ energy levels. The non-dimensionalized parameters for the model, consistent with making the rotating frame approximation, are summarized in Table \ref{table:parameters_JC}. 
\begin{table}[hbt]
\centering 
\caption{Parameters for Jaynes-Cumming Model.}
\label{table:parameters_JC}
\begin{tabular}{|c|c|c|c|c|c|c|c|c|c|}
\hline  
$\omega_0/2\pi$ & $\omega_1/2\pi$ & $\xi_0/2\pi$& $\xi_1/2\pi$& $J_{01}$ & $\xi_{01}$& $T_{00}$&$T_{01}$& $T_{10}$&$T_{11}$\\
\hline  
$0$ & $0$ &  $0$& $0$&1.0&$0$&$\infty$&4500&$\infty$& $\infty$  \\
\hline  
\end{tabular}
\end{table}

\begin{figure}[hbt]
\centering
\includegraphics[width=0.5\textwidth, trim=0.0cm 0cm 0.0cm 0cm, clip=true]{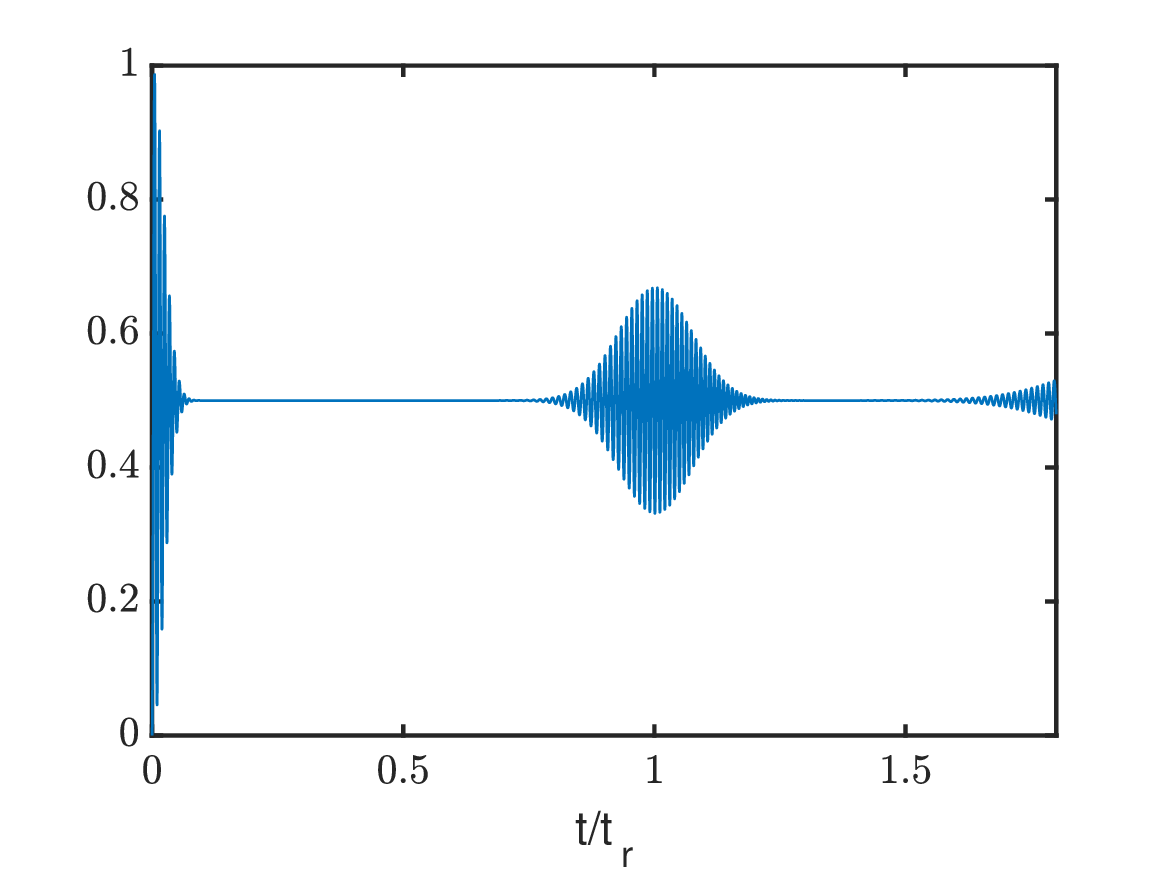}
\caption{Displayed is the population as a function of time for the ground state in the qubit. The revival phenomena is clearly manifested around the revival time $t_r$.
\label{fig:population_JC}}
\end{figure}

Following \cite{PhysRevA.44.5913}, at the initial time the atom is in the excited state and the cavity is in the coherent state 
\[
{\bf v} \sim \sum_{n=0}^{m-1} \frac{|v|^n}{\sqrt{n!}} {\bf e}_n,   
\]
that is, the initial density matrix is $\rho = VV^\dagger$, with 
\[
V =  \begin{pmatrix}
0 \\
1
\end{pmatrix} \otimes \frac{{\bf v}}{\|{\bf v} \|}.
\]
Here we choose $v = \sqrt{m/3}$ so that the last terms in the sum for ${\bf v}$ are small also for moderate $m$. With these choices, without control, and for large $m$ the populations in the atom collapses to $1/2$ but revives at $t_r = 2 \pi |v| / \lambda$, see \cite{PhysRevA.44.5913}. This behavior is displayed in Figure \ref{fig:population_JC}.

\begin{figure}[hbt]  
\centering
\includegraphics[width=0.47\textwidth, trim=0.0cm 0cm 0.0cm 0cm, clip=true]{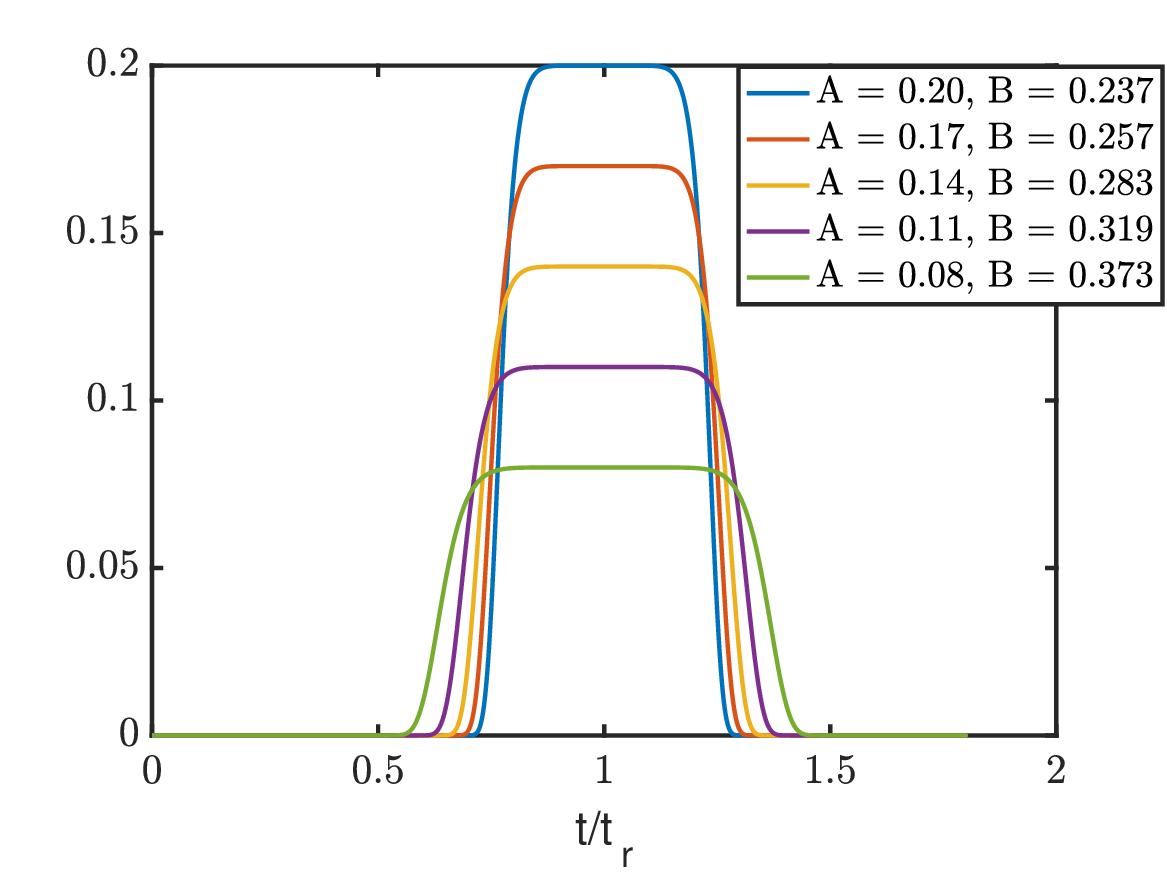}
\includegraphics[width=0.43\textwidth, trim=0.0cm 0cm 1.0cm 1cm, clip=true]{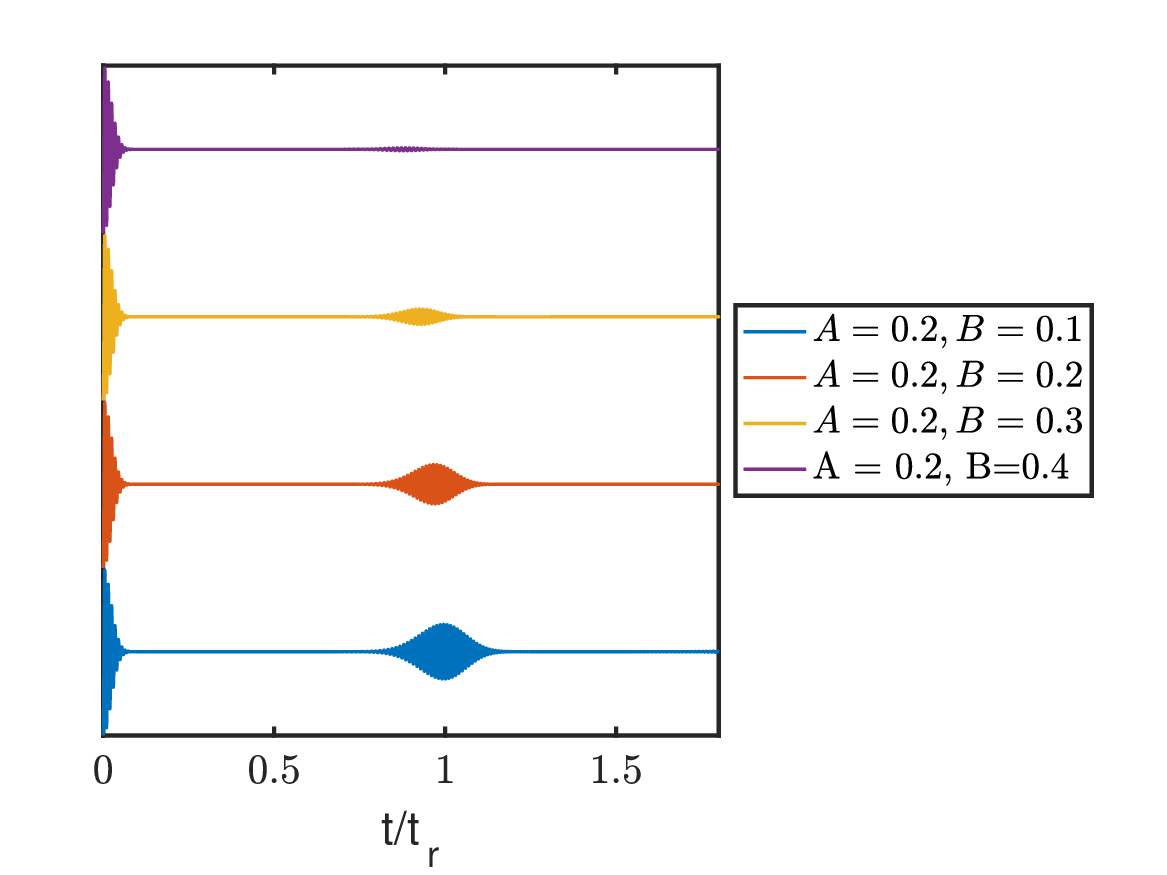}
\caption{Left: the control function (\ref{eq:JCfun}) for different values of $A$ and $B$; Right: the population in the ground state for the qubit as a function of time. Note that the three top curves are offset in the vertical direction by 1,2,3. \label{fig:JC_1A_4Bcontrol}}
\end{figure}

We now investigate if the revival can be suppressed by a suitably chosen control protocol that drives the qubit population to a steady value of $1/2$ and maintains it throughout the revival timescale $t_r$. We choose a control function consisting of the following smooth, rapidly decaying  function
\begin{equation} \label{eq:JCfun}
f_c(t)=Ae^{-\left(\frac{t/t_r-1}{B}\right)^{10}}.
\end{equation}
Here $A,B$ are two parameters. We plot this control for different choices of parameters $A$ and $B$ in Figure \ref{fig:JC_1A_4Bcontrol}.

To quantify how the two parameters $A,B$ affect the revival in a time interval around $t_r$, we define a cost value function 
\[
C(A,B):=\int e^{ -\left(\frac{t/t_r-1}{0.6}\right)^{10}} (P(t)-\frac{1}{2})^2 dt,
\]
where $P(t)$ denotes the population of the qubit in the excited state (this depends implicitly on the parameters $A,B$). 
We estimate $A$ and $B$ by minimizing  the cost function $C(A,B)$. The data set consists of $40\times 60$ evaluations of $C(A,B)$ on a uniform grid, with $A\in\{0.00,0.01,\ldots,0.40\}$ and $B\in\{0.01,0.02,\ldots,0.60\}$. For each parameter pair $(A,B)$, we compute the population dynamics $P(t)$ using the explicit fourth order scheme with $6400$ time steps. 

\begin{figure}[htb] 
\centering
\includegraphics[width=0.5\textwidth, trim=2.0cm 6cm 2.0cm 7cm, clip=true]{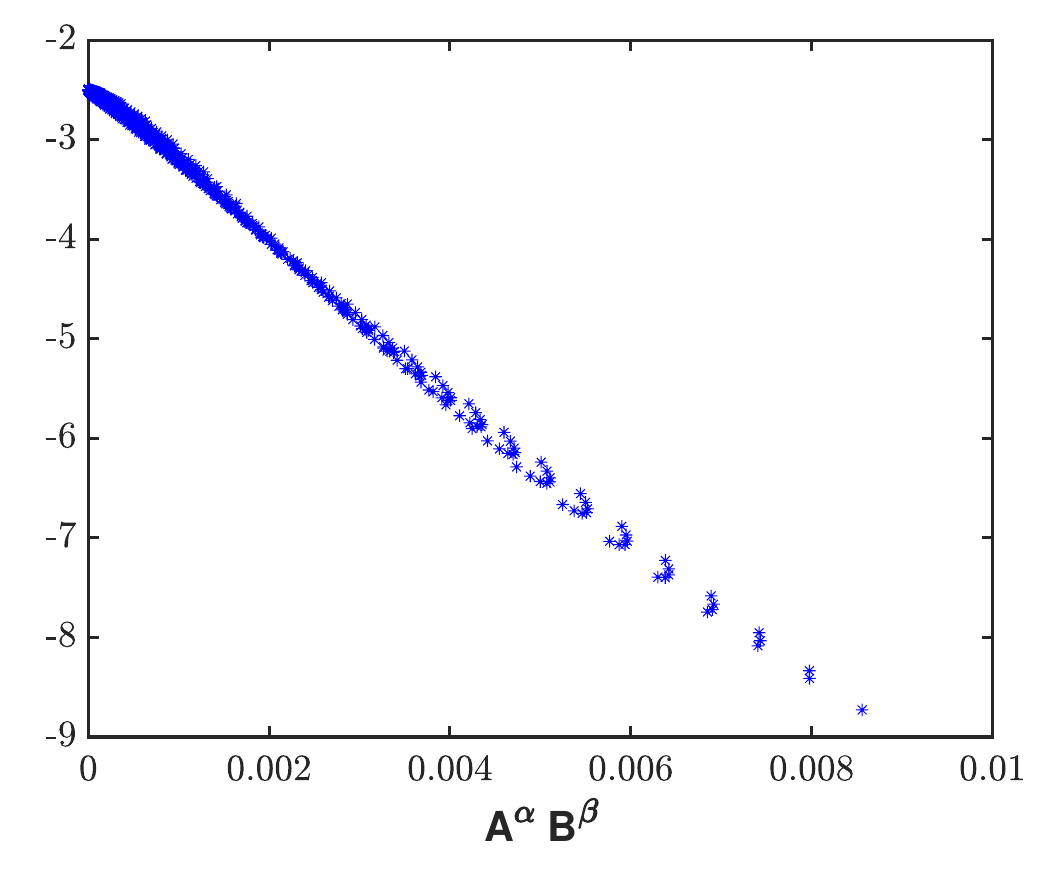}
\caption{Data after fitting with $\alpha=1.3728$ and $\beta=2.7844$, where the  $y$ axis represents $\log(C(A,B))$.
\label{fig:data_fitting_best}}
\end{figure}

Figure~\ref{fig:JC_1A_4Bcontrol} illustrates the population dynamics for four representative $(A,B)$ pairs. Our results demonstrate that increasing $B$ leads to stronger suppression of the revival phenomenon. We hypothesize that the cost depends on $ A $ and $ B $ through the relation
\[
\log(C(A, B)) = \theta_0 + \theta_1 A^\alpha B^\beta,
\]
where $ \theta_0, \theta_1, \alpha, \beta $ are parameters to be estimated from the data via nonlinear least squares (using Matlab's toolbox). This model provides a good fit for the region where $ A \leq 0.2 $ and $ B \leq 0.4 $. The estimated values of $ \alpha, \beta, \theta_0, \theta_1 $, along with their 95\% confidence intervals, are reported in Table~\ref{table:data_fitting_cost_JC}. Using the fitted parameters \( \alpha = 1.3728 \) and \( \beta = 2.7844 \), we demonstrate the model's performance in Figure~\ref{fig:data_fitting_best}. Notably, all \((A,B)\) pairs shown to the left  in Figure~\ref{fig:JC_1A_4Bcontrol}  yield identical values of the composite parameter \(A^\alpha B^\beta =0.002\). Consequently, each control function in Figure~\ref{fig:JC_1A_4Bcontrol} generates approximately equivalent cost values.

\begin{table}[hbt]
\centering 
\caption{Data fitting for the cost value function.}
\label{table:data_fitting_cost_JC}
\begin{tabular}{|c|c|c|}
\hline
& Lower Bound & Upper Bound\\
\hline
$\alpha$ & 1.3623 & 1.3834\\
\hline  
$\beta$ & 2.7647 &  2.8041 \\
\hline  
$\theta_0$ &-2.4757& -2.4643\\
\hline  
$\theta_1$&-792.9967& -745.8285  \\
\hline  
\end{tabular}
\end{table}

The data fitting result  suggests that the cost function is exponentially decaying and the decay rate is affected by the value of $A^\alpha B^\beta$. This phenomenon is clearly demonstrated to the right in Figure~\ref{fig:JC_1A_4Bcontrol}, which compares four parameter pairs with fixed $A=0.2$ and varying $B$ values ($0.1$, $0.2$, $0.3$, $0.4$). The pair $(A,B) = (0.2,0.4)$ yields the largest value of the composite parameter $A^\alpha B^\beta$, which correlates with the strongest suppression of revival phenomena in our simulations.

\section{Conclusions}\label{sec: conclusions}
We develop a family of CPTP schemes for solving the Lindblad equation with a time-dependent Hamiltonian that models the evolution of open quantum systems. By NPI, arbitrary order in time accuracy is achieved. We note that our method by NPI is naturally an embedded integrator that can be used to perform time-adaptive calculations. Furthermore, we carry out stability analysis to quantify the stability conditions for schemes in various order. To address the computational challenges for large systems, we used low-rank truncation by SVDs. We validate the proposed methods through  numerical experiments, demonstrating both high accuracy and computational efficiency in a variety of representative quantum dynamical scenarios.  One natural extension of this work is to consider larger systems using tensor network representations for the columns of $V$. We are currently pursuing this research avenue. Another extension is to use that the method is embedded to design a time step adaptive method.

\bibliographystyle{plain}
\bibliography{ref}

@book{wanner1996solving,
  title={Solving ordinary differential equations II},
  author={Wanner, Gerhard and Hairer, Ernst},
  volume={375},
  year={1996},
  publisher={Springer Berlin Heidelberg New York}
}

@article{PhysRevA.44.5913,
  title = {Atom- and field-state evolution in the {J}aynes-{C}ummings model for large initial fields},
  author = {Gea-Banacloche, Julio},
  journal = {Phys. Rev. A},
  volume = 44,
  issue = 9,
  pages = {5913--5931},
  numpages = 0,
  year = 1991,
  month = {Nov},
  publisher = {American Physical Society}
}

@article{HOHO,
      title={High-Order {H}ermite {O}ptimization: {F}ast and {E}xact {G}radient {C}omputation in {O}pen-{L}oop {Q}uantum {O}ptimal {C}ontrol using a {D}iscrete {A}djoint {A}pproach}, 
      author={Spencer Lee and Daniel Appel\"{o}},
      year={2025},
      journal={Journal of Computational Physics (submitted) arXiv 2505.09857}
}

@article{robin2025unconditionally,
  title={Unconditionally stable time discretization of {L}indblad master equations in infinite dimension using quantum channels},
  author={Robin, R{\'e}mi and Rouchon, Pierre and Sellem, Lev-Arcady},
  journal={arXiv preprint arXiv:2503.01712},
  year={2025}
}

@article{puzynin2000magnus,
  title={Magnus-factorized method for numerical solving the time-dependent {S}chr{\"o}dinger equation},
  author={Puzynin, IV and Selin, AV and Vinitsky, SI},
  journal={Computer physics communications},
  volume={126},
  number={1-2},
  pages={158--161},
  year={2000},
  publisher={Elsevier}
}

@article{dutt2000spectral,
  title={Spectral deferred correction methods for ordinary differential equations},
  author={Dutt, Alok and Greengard, Leslie and Rokhlin, Vladimir},
  journal={BIT Numerical Mathematics},
  volume={40},
  number={2},
  pages={241--266},
  year={2000},
  publisher={Springer}
}

@article{schratz2021low,
  title={Low-regularity integrators for nonlinear {D}irac equations},
  author={Schratz, Katharina and Wang, Yan and Zhao, Xiaofei},
  journal={Mathematics of Computation},
  volume={90},
  number={327},
  pages={189--214},
  year={2021}
}

@article{baumstark2018uniformly,
  title={Uniformly accurate exponential-type integrators for {K}lein-{G}ordon equations with asymptotic convergence to the classical {NLS} splitting},
  author={Baumstark, Simon and Faou, Erwan and Schratz, Katharina},
  journal={Mathematics of Computation},
  volume={87},
  number={311},
  pages={1227--1254},
  year={2018}
}

@article{cai2019uniformly,
  title={Uniformly accurate nested {P}icard iterative integrators for the {D}irac equation in the nonrelativistic limit regime},
  author={Cai, Yongyong and Wang, Yan},
  journal={SIAM Journal on Numerical Analysis},
  volume={57},
  number={4},
  pages={1602--1624},
  year={2019},
  publisher={SIAM}
}

@article{chen2025full,
  title={Full{-} and low-rank exponential midpoint schemes for forward and adjoint {L}indblad equations},
  author={Chen, Hao and Borzi, Alfio},
  journal={arXiv preprint arXiv:2506.00346},
  year={2025}
}

@article{lawson1967generalized,
  title={Generalized {R}unge-{K}utta processes for stable systems with large {L}ipschitz constants},
  author={Lawson, J Douglas},
  journal={SIAM Journal on Numerical Analysis},
  volume={4},
  number={3},
  pages={372--380},
  year={1967},
  publisher={SIAM}
}

@article{bachmayr2025iterative,
  title={Iterative thresholding low-rank time integration},
  author={Bachmayr, Markus and Dolbeault, Matthieu and Sachsenmaier, Polina},
  journal={arXiv preprint arXiv:2507.15848},
  year={2025}
}

@article{li2024high,
  title={{H}igh-{O}rder {I}mplicit {L}ow-{R}ank {M}ethod with {S}pectral {D}eferred {C}orrection for {M}atrix {D}ifferential {E}quations},
  author={Li, Shun and Jiang, Yan and Cheng, Yingda},
  journal={arXiv preprint arXiv:2412.09400},
  year={2024}
}

@article{appelo2024kraus,
      title={Kraus is {K}ing: High-order {C}ompletely {P}ositive and {T}race {P}reserving {(CPTP)} {L}ow {R}ank {M}ethod for the {L}indblad {M}aster {E}quation}, 
      author={Daniel Appel\"{o} and Yingda Cheng},
      journal={Journal of Computational Physics},
volume = {534},
pages = {114036},
year = {2025}}

@article{manzano2020short,
  title={A short introduction to the {L}indblad master equation},
  author={Manzano, Daniel},
  journal={Aip advances},
  volume={10},
  number={2},
  year={2020},
  publisher={AIP Publishing}
}

@article{cao2025structure,
  title={Structure-preserving numerical schemes for {L}indblad equations},
  author={Cao, Yu and Lu, Jianfeng},
  journal={Journal of Scientific Computing},
  volume={102},
  number={1},
  pages={1--34},
  year={2025},
  publisher={Springer}
}

@article{gunther2021quantum,
  title={Quantum optimal control for pure-state preparation using one initial state},
  author={G{\"u}nther, Stefanie and Petersson, N Anders and DuBois, Jonathan L},
  journal={AVS Quantum Science},
  volume={3},
  number={4},
  year={2021},
  publisher={AIP Publishing}
}

@article{steinbach1995high,
  title={High-order unraveling of master equations for dissipative evolution},
  author={Steinbach, J and Garraway, BM and Knight, PL},
  journal={Physical Review A},
  volume={51},
  number={4},
  pages={3302},
  year={1995},
  publisher={APS}
}

@article{jaynes1963comparison,
  title={Comparison of quantum and semiclassical radiation theories with application to the beam maser},
  author={Jaynes, Edwin T and Cummings, Frederick W},
  journal={Proceedings of the IEEE},
  volume={51},
  number={1},
  pages={89--109},
  year={1963},
  publisher={IEEE}
}

@article{lindblad1976generators,
  title={On the generators of quantum dynamical semigroups},
  author={Lindblad, Goran},
  journal={Communications in mathematical physics},
  volume={48},
  pages={119--130},
  year={1976},
  publisher={Springer}
}

@article{chen2024full,
  title={Full{-}  and low-rank exponential {E}uler integrators for the {L}indblad equation},
  author={Chen, Hao and Borz{\`\i}, Alfio and Jankovi{\'c}, Denis and Hartmann, Jean-Gabriel and Hervieux, Paul-Antoine},
  journal={arXiv preprint arXiv:2408.13601},
  year={2024}
}

@article{PhysRevA.76.042319,
  title = {Charge-insensitive qubit design derived from the {C}ooper pair box},
  author = {Koch, Jens and Yu, Terri M. and Gambetta, Jay and Houck, A. A. and Schuster, D. I. and Majer, J. and Blais, Alexandre and Devoret, M. H. and Girvin, S. M. and Schoelkopf, R. J.},
  journal = {Phys. Rev. A},
  volume = {76},
  issue = {4},
  pages = {042319},
  numpages = {19},
  year = {2007},
  month = {Oct},
  publisher = {American Physical Society},
  doi = {10.1103/PhysRevA.76.042319},
  url = {https://link.aps.org/doi/10.1103/PhysRevA.76.042319}
}

@article{PhysRevLett.101.080502,
  title = {Controlling the {S}pontaneous {E}mission of a {S}uperconducting {T}ransmon {Q}ubit},
  author = {Houck, A. A. and Schreier, J. A. and Johnson, B. R. and Chow, J. M. and Koch, Jens and Gambetta, J. M. and Schuster, D. I. and Frunzio, L. and Devoret, M. H. and Girvin, S. M. and Schoelkopf, R. J.},
  journal = {Phys. Rev. Lett.},
  volume = {101},
  issue = {8},
  pages = {080502},
  numpages = {4},
  year = {2008},
  month = {Aug},
  publisher = {American Physical Society},
  doi = {10.1103/PhysRevLett.101.080502},
  url = {https://link.aps.org/doi/10.1103/PhysRevLett.101.080502}
}

@article{
doi:10.1126/science.1231930,
author = {M. H. Devoret  and R. J. Schoelkopf },
title = {Superconducting {C}ircuits for {Q}uantum {I}nformation: {A}n {O}utlook},
journal = {Science},
volume = {339},
number = {6124},
pages = {1169-1174},
year = {2013},
doi = {10.1126/science.1231930},
URL = {https://www.science.org/doi/abs/10.1126/science.1231930},
eprint = {https://www.science.org/doi/pdf/10.1126/science.1231930}}

\end{document}